\documentclass{amsart}

\usepackage{graphicx}
\usepackage{psfrag}
\usepackage{amsthm}
\usepackage{amscd}
\usepackage{amsmath}
\usepackage{amsfonts}
\usepackage{amstext}

\begin{document}

\title{Renormalization theory for multimodal maps} 
\author{Daniel Smania }
\footnote{ This work has been partially supported by CNPq and FAPERJ grant E-26/151/462/99.}
\address{Instituto de Matem\'atica Pura e Aplicada, 
    Estrada Dona Castorina, 110, 
               Jardim Bot\^anico,
               CEP 22460-320, 
               Rio de Janeiro-RJ,
               Brazil. }
\email{smania@impa.br}

\begin{abstract} We study the dynamics of the renormalization operator for multimodal maps. In particular, we prove the exponential convergence of this operator for infinitely renormalizable maps with same bounded combinatorial type. \end{abstract}

\maketitle

\newcommand{\co}{\mathbb{C}}
\newcommand{\incl}[1]{i_{U_{#1}-Q_{#1},V_{#1}-P_{#1}}}
\newcommand{\inclu}[1]{i_{V_{#1}-P_{#1},\co-P}}
\newcommand{\func}[3]{#1\colon #2 \rightarrow #3}
\newcommand{\norm}[1]{\left\lVert#1\right\rVert}
\newcommand{\norma}[2]{\left\lVert#1\right\rVert_{#2}}
\newcommand{\hiper}[3]{\left \lVert#1\right\rVert_{#2,#3}}
\newcommand{\hip}[2]{\left \lVert#1\right\rVert_{U_{#2} - Q_{#2},V_{#2} - P_{#2}}}
\newtheorem{lem}{Lemma}[section]
\newtheorem{defi}{Definition}[section]
\newtheorem{rem}{Remark}[section]
\newtheorem{mle}{Main Lemma}[section]
\newtheorem{thm}{Theorem}
\newtheorem{cor}{Corollary}[section]
\newtheorem{pro}{Proposition}[section]
\newtheorem{conj}{Conjecture}

\section{Introduction}

The renormalization theory has a long history: beginning with observable universality properties  and conjectural explanation of these observations in families of unimodal maps, by Feigenbaum and Collet-Tresser. O. Lanford proposed the existence of a hyperbolic horseshoe to the renormalization operator. Similar conjectures was done for critical circle maps and for bimodal maps (\cite{Mackay}).

A new step in the renormalization theory was attained by Sullivan's work(\cite{S}): new tools was introduced, like quasiconformal deformation methods and a fruitful analogy with the theory of Kleinian groups. McMullen(\cite{Mc2}) proved the exponential convergence of the renormalization operator and Lyubich (\cite{Lyu99}) its hyperbolicity (in the space of quadratic-like maps). 

Our intend is construct the foundation of the renormalization theory for multimodal maps. The pioneer in this issue was J. Hu (see \cite{H95} and  \cite{HU98}), which studied the renormalization operator (compactness and convergence) for bimodal maps in the Epstein class. We will study the convergence of the  renormalization operator using the methods introduced by the cited authors for unimodal maps. 

\subsection{Outline of paper}  

In the section \ref{multimodal} we will introduce the most important object in the paper: multimodal maps of type $n$. These maps are maps obted of compositions of unimodal maps. Indeed, deep renormalizations of multimodal maps are multimodal maps of type $n$, so there are not loss of generality in restrict our study for these maps. Furthermore, these maps have a nice structure: in particular, we can define the combinatorial type of a renormalization, give explicit rules to the compositions of combinatorial types and realize any combinatorial type in sufficiently rich families. This is done in section \ref{comb}. In section \ref{pol} we study polynomials which are compositions of quadratic polynomials. We prove that the locus of connectivity of this family is compact and give a criteria to decide when a polynomials is a composition of quadratic ones. In section \ref{pol-like} we introduce the polynomial like maps of type $n$, and we prove that these maps are hybrid conjugated with compositions of quadratic polynomials. We study also compact subsets in the space of polynomial like maps. In section \ref{ren} we define the complex version of renormalization and prove the 'small Julia set everywhere' theorem, which implies, in particular, that infinitely renormalizable polynomial like maps of type $n$ with bounded combinatorics does not support non trivial Beltrami fields in its Julia set. This result will be used in section \ref{hyb}, where we prove that infinitely renormalizable real polynomials of type $n$  with same bounded combinatorial type are hybrid conjugated. As a corollary, we obtain that the  set of  infinitely renormalizable real polynomials of type $n$  with combinatorial type bounded by a constant $C$ is a Cantor set. In the section \ref{tow} we define the McMullen's towers and prove it rigidity. The theory is quite similar to the unimodal case and it implies the convergence of renormalization. Finally, in the section \ref{expo} we prove, using the McMullen's theory of dynamic inflexibility, the exponential convergence of the renormalization operator. In the apendice we collect some results about fixed-point theory,  a special kind of Riemann surface and quasiconformal theory.

\subsection{Multimodal maps}\label{multimodal}

A multimodal map $f\colon I \rightarrow I$, $I=[-1,1]$, is a smooth map with a finite number of critical points, all of them local maximum or local minimum, and such that $f(\partial I) \subset \partial I$. Here we will be interested in more specific kinds of multimodal maps:

\begin{defi}\label{multin} We say that $f$ is a \textbf{multimodal map of type $n$} if it can be written as a composition of $n$ unimodal maps: to be more precise, there exist  maps $f_1, \dots, f_n$ with the following properties
\begin{enumerate}
\item $f_i\colon I \rightarrow I$ has an unique critical point (a maximum) and  $f_i(\partial I) \subset \partial I$.
\item $f = f_n \circ \dots \circ f_1$.
\item If $c_i$ is the critical point of $f_i$, then $f_i(c_i) \geq c_{i+1 \text{ mod }n}$.
\end{enumerate}
The $n$-uple $(f_1,...,f_n)$ is a \textbf{decomposition} of $f$. Clearly $f$ has many decompositions. For each decomposition of $f$ we can associate an \textbf{extended map} $F$ defined on $I_n = \{(x,i) \colon x \in I, 1 \leq i \leq n  \}$ (in other words, $I_n$ is a disjoint union of $n$ copies of $I$) by
\begin{equation}
F(x,i) = (f_i(x), i + 1 \text{ mod }n)
\end{equation}
\end{defi}

For $(x,i), (y,j) \in I_n$, we say that $(x,i) < (y,j)$ if $i = j$ and $x < y$. The  intervals of  $I_n$ are the sets $J \times \{i\}$, for some $ J \subset I$ and $i < n$. If $c_i$ is the critical point of $f_i$, denote $C(F)=\{(i,c_i)  \}_i$.

In \cite{Sm}, we proved that deep renormalizations in infinitely renormalizable multimodal maps are multimodal maps of type $n$. This is the reason to restrict our attention for this kind of map. 
\begin{defi} We say that $J$ is a \textbf{$\mathbf{k}$-periodic interval} to the extended map  $F$ if
\begin{itemize}
\item $(c_1,0) \in J$ ($c_i$ is the critical point of $f_i$),
\item $\{ J, F(J), \dots, F^{k-1}(J) \}$ is an union of intervals with disjoint interior,
\item The union of intervals in the above family contains $\{ (c_i,i) \}$,
\item $F^k(J) \subset J$, for $k > n$. 
\end{itemize}
We will call $k$ the \textbf{period} of $J$.
\end{defi}

Suppose that there exists a $k$-periodic interval for $F$. Let  $P$ be the maximal interval which is a $k$-periodic interval for $F$. Then $F^{k}(\partial P) \subset \partial P$. We say that $P$ is a \textbf{restrictive interval} for $F$ of period $k$. Note that if $P$ and $\tilde{P}$ are, respectively, restrictive intervals for $F$ of period $k$ and $\tilde{k}$, $k < \tilde{k}$, then $\tilde{P} \subset P$.  Let $P$ be a restrictive interval and  $0 = \ell_1 < \dots < \ell_n$ the times such that $(c_i,i) \in F^{\ell_j}(P)$ for some $i$. Let $P_j$ be the symmetrization of $F^{\ell_j}(P)$ in relation to $(c_i,i)$. Observe that $P_j$ contains a periodic point in its boundary. Let $A_{P_j}$ be the affine map which maps $P_j$ to $I$ and this periodic point to $-1$. Then  $g_j =  A_{P_{j+1}} \circ F^{\ell_{j+1}-\ell_{j}}\circ A_{P_j}^{-1}$ is a unimodal map. Hence $g = A_{P_1} \circ F^{k} \circ A_{P_1}^{-1}$ is a multimodal map of type $n$ with decomposition $(g_1,\dots,g_n)$. If $k > n$ is the minimal number such that $F$  admits a restrictive interval of period $k$, the map $g$ is called the \textbf{renormalization} of $f$, and denoted by $R(f)$. Indeed, it is easy observe that the definition of $R(f)$ does not depend on the decomposition.

The map $R(f)$ can be renormalizable again and so on. If this process never finished, we say that $f$ is \textbf{infinitely renormalizable}. Denote by $P^k_0$ the restrictive interval associate to the $k$-th renormalization $R^k(f)$. If $q \in C(F)$, denote by the corresponding capital letter $Q_0^k$ the symmetrization of the interval $F^{\ell}(P_0^k)$ which contains $q$. We reserve the letter $p$ for $(c_1,1)$. The critical point $r$ for $F$ will be the \textbf{successor} of the critical point $q$ at level $k$ if $r \in F^\ell(Q_0^k)$, for the minimal $\ell$ so that $F^\ell(Q_0^k)$ contains a critical point. Define $n^k_r = \ell$. Then, for any $q \in C(F)$, $k \in \mathbb{N}$ and $i < n^k_r$, there exists an interval $R^k_{-i}$ so that
\begin{itemize}
\item $F^i$ is monotone in $R^k_{-i}$,

\item $F^i(R^k_{-i})=R^k_{0}$,
 
\item The interval $F^{n^k_r - i}(Q^k_0)$ is contained in $R^k_{-i}$.  
\end{itemize}
For details, see \cite{Sm}.

Denote by $N_k$ the period of the restrictive interval $P^k_0$. We say that $f$ has \textbf{C-bounded combinatorics} if  $N_{k+1}/N_k \leq C$.   
 
\section{Combinatorial results}\label{comb}

\begin{defi} Let $f$ be a multimodal map of type $n$. Let $(f_1,...,f_n)$ be a decomposition of $f$. Let  $x$ be a point in the domain of the extended map $F$ associate with this decomposition. The \textbf{itinerary} of $x$ with respect to the decomposition $(f_1,...,f_n)$ is the infinity word $\ell_0(x) \ell_1(x) \dots \ell_i(x) \dots$, with $\ell_i(x) = L, C, R$ satisfying
\begin{equation}
\ell_i(x) = \begin{cases} R&    \text{ if } F^i(x) > c_j, \text{ for some j;}  \\
                       C&     \text{ if } F^i(x) = c_j, \text{ for some j;} \\
                       L&     \text{ if }  F^i(x) < c_j, \text{ for some j.} \end{cases}
\end{equation}
Let $(x,i)$ be  a point of $I_n$. The \textbf{ inner itinerary} of $(x,i)$ is the finite word $\ell_0(x,i) \ell_2(x,i) \dots \ell_{n - i}(x,i)$. 
\end{defi}

Let $f$ be a multimodal map of type $n$ with $n$ distinct critical values. Order the critical points of $f$ , $a_1 < \dots a_{k_f}$, $k_f < 2^n$,  and let $v_1, \dots, v_n$ be the critical values of $f$. We associate to $f$ the structure $<k_f,\psi_f>$, where $\psi_f$ is the map of $\{ i \in \mathbb{N}\colon 1 \leq i \leq k_f \}$ into $\{ i \in \mathbb{N}\colon 1 \leq i \leq n \}$ such that $\psi_f(i)=j$ iff $f(a_i)=v_j$.

The definition of the  structure $<k_f,\psi_f>$ does not depend on the decompositions, but the inner itinerary of the critical point of a extended map  does. Indeed, as it is natural, these two combinatorial informations are completely equivalents. The proof is quite boring to read (or write), but elementary.

\begin{lem}\label{inner} Let $f$ be a multimodal map of type $n$ with $n$ distinct critical values. Then the inner itinerary of the critical points of the extended map $F$ depends only on the structure $<k_f,\psi_f>$. In particular, the inner itinerary does not depend on the decomposition of $f$. Furthermore, $<k_f,\psi_f>$ can be determinated by the inner itinerary associate to a decomposition.
\end{lem}

\begin{proof} 
first part: Let $(f_1,...,f_n)$ be an arbitrary decomposition of $f$, $c_i$ the critical point of $f_i$ and let $a_1 < \dots < a_k$, $k < 2^n$ be the critical points of $f$. Observe that for each $a_i$\, there is an unique letter $C$ in its itinerary, since $f$ has $n$ critical values. Hence $a_1$ has itinerary $LL\dots LLC$ (by item 3 in definition \ref{multin}). In particular, $f(a_1)=f_{n}(c_n)$. Suppose by induction that we know the inner itinerary of $a_i$, for $i \leq j$. Let $\tilde{\ell}_1\dots \tilde{\ell}_n$ be the itinerary of $a_j$ and $\ell_1 \dots \ell_n$ the itinerary of $a_{j+1}$. The itinerary of $a_{j+1}$ is equal to the itinerary of  $a_{j}$, except in two positions: if $\tilde{\ell}_r = C$ then clearly $\ell_r \neq C$ and if $\ell_r = C$ then $\tilde{\ell_r} \neq C$. To find the position of the letter $C$, there are two cases: if there is $i \leq j$ such that  $f(a_{j+1})=f(a_i)$, we conclude that the letter $C$ happens in the itinerary of $a_{j+1}$ in the same position that in the inner itinerary of $a_{i}$.  If $f(a_{j+1})\neq a_i$ for $i \leq j$, then the letter $C$ happens in the itinerary of $a_{j+1}$ exactly in the position $min\{ r \colon C \text{ happens in the r-th position in the itinerary of some } a_i, \text{ }i \leq j  \} - 1$. Now the letter in the $r$-th position in the itinerary of $a_{j+1}$, where $\tilde{\ell}_r=C$, will be $R$ if $j+1$ is even and the word $\tilde{\ell}_1 \dots \tilde{\ell}_{r-1}R\tilde{\ell}_{r+1}\dots \tilde{\ell}_n$ has  an even number of letters $R$; or  $j+1$ is odd and the word $\tilde{\ell}_1 \dots \tilde{\ell}_{r-1}R\tilde{\ell}_{r+1}\dots \tilde{\ell}_n$ has  an odd number of letters $R$. Otherwise $\ell_r = L$.

second part: We are going to prove that $<k_f,\psi_f>$ can be deduced of the inner itinerary. We will get a proof by induction in $n$. Suppose that for any multimodal map $g$ of type $n$ and a decomposition $(g_1,\dots,g_n)$, it is possible obtain $<k_g,\psi_g>$ of the inner itinerary of the critical points of the associate extended map $G$. Consider a multimodal map $f$ of type $n+1$ with decomposition $(f_1,\dots,f_{n+1})$. Then the structure $<k_g,\psi_g>$ associate to the map $g= f_{n+1} \circ \dots \circ f_2$ can be deduced of the inner itinerary of $f$, by induction hypothesis. Let $a^g_1 < \dots < a^g_{k_g}$ be the critical points of $g$. By the proof of the first part, the inner itinerary of the points in each interval between the critical points of $g$ can be deduced of the structure $<k_g,\psi_g>$. Select the interval $J=[a^g_j,a^g_{j+1}]$ such that any point in $J$ has the same inner itinerary that $f_1(0)$ with respect to $g$. Define $C = \{i \in \mathbb{Z} \colon  |i| \leq j \}$ and the map $\psi \colon C \rightarrow \{i \colon 1 \leq i \leq n+1 \}$ by $\psi(i)= \psi(|i|)$, if $i \neq 0$ and  $\psi(0)= n+1$. Then $k_f = \# C$ and $\psi_f(i)=\psi(i + k_f + 1)$.

\end{proof}
\begin{rem} We will not give details, but  $<k_f,\psi_f,\phi_f>$, where $\phi_f(i)= r$ if the $i$-th critical point of $f$ has criticality $r$, is a combinatorial information equivalent to the inner itinerary of the critical points of a extended map for $f$, even if $f$ does not have $n$ critical values.  \end{rem}

In an analogous way,  the itinerary of an interval $J \subset I_n$  is defined by
\begin{equation}
\ell_i(J) = \begin{cases} R&    \text{ if } F^i(J) > c_j, \text{ for some j;}  \\
                       C&     \text{ if } F^i(J) \text{ contains } c_j, \text{ for some j;} \\
                       L&     \text{ if }  F^i(J) < c_j, \text{ for some j.} \end{cases}
\end{equation} 
The inner itinerary of $J \subset I \times \{ i\}$ is the finite word $\ell_0(J) \dots \ell_{n - i}(J)$.

\begin{cor} Let $J$ be an interval in $I$, and $f$ a multimodal map of type $n$. Then the inner itinerary of $J$ does not depend on the decomposition.\end{cor}

By the above corollary, if $J$ is a $k$-periodic interval  for some decomposition of $f$ then $J$ is $k$-periodic for all decomposition. In particular the maximal interval $P^1_0$ does not depend on the decomposition and so do the renormalization of $f$. Again by the previous theorem, the order of the intervals in the orbit of $J$ by an extended map $F$ in the $n$ copies on $I$ does not depend on the decomposition. 

\begin{cor} The itinerary of the point in $I \subset I_n$ with respect an extended map $F$ does not depend on the decomposition.
\end{cor}

The \textbf{signal} of a finite word $\omega = \ell_0\dots \ell_k$, $sgn(\omega)$, will be $1$, if there exists an even number of letters $R$ in $\omega$, or $-1$ otherwise. We will apply the signal function only on words which do not contain the letter $C$ (then we say that the word $\omega$ is \textbf{pure}). Provide the set of finite pure words  with the following order $\prec$, defined by
\begin{itemize} 
\item Provide the set of words with length one with the order: $L \prec C \prec R$.
\item If $\omega = \ell_0 \dots \ell_n$ and $\tilde{\omega} = \tilde{\ell}_0 \dots  \tilde{\ell}_n$ are such that  $\ell_0 \dots \ell_j = \tilde{\ell}_0 \dots  \tilde{\ell}_j$, but $\ell_{j+1} \neq \tilde{\ell}_{j+1}$, then
\begin{itemize}
 \item $\omega \prec \tilde{\omega}$ if $sgn(\ell_0 \dots \ell_j)= 1$ and $\ell_{j+1} = L$; or $sgn(\ell_0 \dots \ell_j)= -1$ and $\ell_{j+1} = R$.
 \item Otherwise $\tilde{\omega} \prec \omega$.
\end{itemize}
\end{itemize}

This is the usual order to words with two symbols which ocurrs in the study of unimodal maps (see, e.g, \cite{MS}).

\begin{lem}\label{preorder} Let $x, y \in I$. Assume that the pure itineraries $\ell_0(x) \dots \ell_j(x)$ and $\ell_0(y) \dots \ell_j(y)$ are distinct. Then $x < y$ if $\ell_0(x) \dots \ell_j(x) \prec  \ell_0(y) \dots \ell_j(y)$.
\end{lem}
\begin{proof} The proof is easy.
\end{proof}

Let $f$ and $\tilde{f}$ be multimodal maps of type $n$ with decompositions $(f_1,\dots,f_n)$ and $(\tilde{f}_1,\dots,\tilde{f}_n)$. If $c_i$ (resp. $\tilde{c}_i$) is  the critical point of $f_i$ (resp. $\tilde{f}_i$), define $v_i = f_n \circ \dots \circ f_i (c_i)$ (resp. $\tilde{v}_i = \tilde{f}_n \circ \dots \circ \tilde{f}_i (\tilde{c}_i)$. 

\begin{lem}\label{pos} Let $f$ and $\tilde{f}$ be multimodal maps of type $n$ with the same critical itinerary and such that $v_i < v_j$ iff $\tilde{v}_i < \tilde{v}_j$. Let $H_0, H_1 \colon I \rightarrow I$ be increasing continuous functions such that $H_1 \circ f = \tilde{f} \circ H_0$. Then, for  $y \in I$:
\begin{enumerate}
\item $H_0(f^{-1}(y))= \tilde{f}^{-1}(H_1(y))$,
\item For each pure word $\ell_0 \cdots \ell_{n-1}$ there is at most one point $x \in f^{-1}(y)$ such that $\ell_i(x)=\ell_i$.
\item There is a point $x \in f^{-1}(y)$ with inner itinerary $\ell_0 \cdots \ell_{n-1}$  iff there is a point $\tilde{x} \in \tilde{f}^{-1}(H_1(y))$ with the same itinerary. Furthermore $H_0(x)=\tilde{x}$.
\end{enumerate}
\end{lem}
\begin{proof} The item 1 is obvious. To prove 2, let $a_1 <\dots < a_k$ be the critical points of $f$ and $\tilde{a}_1< \dots <\tilde{a}_k$ be the critical points of $\tilde{f}$. Notice that $H_0(a_i)=\tilde{a}_i$ and $H_1(v_j)=\tilde{v}_j$. We saw in the proof of lemma  \ref{inner} that the inner itineraries of $c_i$ and $\tilde{c}_i$ are the same. Thus  one gets that $f(a_i)=v_j$ iff $\tilde{f}(\tilde{a}_i)=\tilde{v}_j$. If the interval $[a_i,a_{i+1}]$ contains a preimage of $y$ then $v_j=f(a_i) \leq y \leq f(a_{i+1})=v_k$. But this occurs iff  $\tilde{v}_j \leq H_1(y) \leq \tilde{v}_k$ and so the interval $[\tilde{a}_i,\tilde{a}_{i+1}]$ contains a preimage of $H_1(y)$. Since the points in $(a_i, a_{i+1})$ and $(\tilde{a}_i, \tilde{a}_{i+1})$ have the same inner itinerary, the proof is finished. 

\end{proof}

\begin{figure}
\centering
\includegraphics[width=0.50\textwidth]{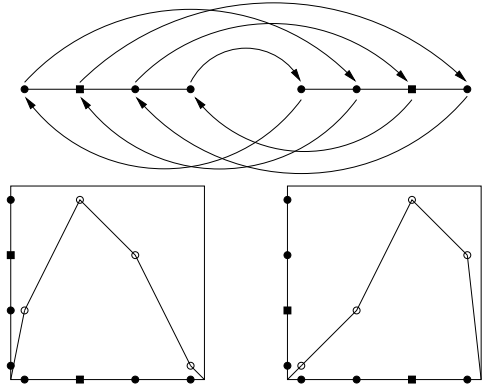}
\caption{In the upper part of the figure we represented a m.c.d.: the elements of $A$ are represented by small discs and squares: the squares are the critical elements. Two elements of $A$ are comparable if they are in the same segment. In this case they respect the order in the real line. The elements are permuted as indicated by the arrows. In the bottom part we see that this m.c.d. can be realizated by a multimodal map of type two.}
\end{figure}

\begin{defi} Denote by  $<A,\prec,A^c, \pi>$ the \textbf{combinatorial data (c.d.)} which contains
\begin{itemize}
\item A finite set $A$ with  $A^c \subset A$. The set $A^c$ is the set of 'critical points' of $A$.

\item   $\prec$ is a transitive, anti-reflexive and anti-symmetric relation under $A$ such that any point in $A$ is comparable with an unique point in $A^c$. Furthermore, the relation 'a is comparable with b on respect to $\prec$', $a \sim b$, is an equivalence relation. We will denote $[x]$ the equivalence classes with respect to $\sim$.  

\item $\pi\colon A \rightarrow A$ is a map with the following property: let $c$ be a critical point in $A^c$ then
\begin{itemize}
\item $a \prec b \prec c$ implies $\pi(a) \prec \pi(b) \prec \pi(c)$.
\item $c \prec b \prec a$ implies $\pi(a) \prec \pi(b) \prec \pi(c)$.
\end{itemize} 
\item For any $a \in A$ there exists $c \in A^c$ so that $\pi^{i}(c)=a$, for some $i \geq 0$.
\end{itemize}

A \textbf{marked combinatorial data (m.c.d.)} will be $<A,\prec,A^c, \pi, c>$ with $c \in A^c$. Two m.c.d. $<\sigma,c>, <\tilde{\sigma}, \tilde{c}>$ will be identified up to bijections $\phi  \colon A \rightarrow \tilde{A}$ satisfying  \begin{enumerate}
\item $\phi(A^c)=\tilde{A}^c$;
\item For $x, y \in A$, $x \prec y$ iff $\phi(x) \tilde{\prec} \phi(y)$;
\item $\pi =  \phi^{-1} \circ \tilde{\pi} \circ \phi$ and
\item $\phi(c)=\tilde{c}$.
\end{enumerate} 
Such map $\phi$ is called an \textbf{isomorphism} between m.c.d. Note that if $\sigma$ and $\tilde{\sigma}$ are two m.c.d. and $\pi$ is \textbf{transitive} (if  $x, y \in A$ then  $\pi^i(x)=y$, for some $i$) then there is at most one isomorphism between $\sigma$ and $\tilde{\sigma}$. A m.c.d is \textbf{essential} if for any pair $a \prec b$ there is $i \geq 0$ so that $\pi^i(a) \preceq d \preceq \pi^i(b)$, for some $d \in A^c$. Clearly if any point in $A^c$ is periodic then  $\sigma$ is essential. \end{defi}

\begin{rem} We can associate two maps to a m.c.d $\sigma$: the \textbf{first entry map to the critical set } $\Pi \colon A \rightarrow A^c$ defined by $\Pi(x)= \pi^i(x)$, where $i \geq 0$ is minimal such that $\pi^i(x) \in A^c$; and the \textbf{first return to the critical set}, defined by $\Pi \circ \pi$.
\end{rem}

\begin{rem} Let $\sigma=<A,\prec,A^c,\pi>$ be a m.c.d. and $x \in A$. The itinerary of $x$ will be the periodic word $\ell_0 \ell_1 \dots \ell_n \dots$ satisfying
\begin{equation} \ell_i = \begin{cases} R&   \text{ if } \pi^i(x) \prec c, \text{ for some $c \in A^c$;} \\
C& \text{ if } c \in A^c; \\
L& \text{ if } c \prec \pi^i(x), \text{ for some $c \in A^c$.} 
\end{cases}
\end{equation}
\end{rem}

\begin{rem} \label{defcf} Let $f$ be a multimodal map of type $n$ and consider $F$ be an extended map induced by a decomposition $(f_1,\dots,f_n)$ of $F$. If $f$ is critically finite, then we can associate to $f$ the following m.c.d.
\begin{itemize}

\item $A^c = \{ (c_i,i) \colon c_i \text{ is the critical point of $f_i$} \}$ and $c = c_1$.

\item $A = \cup_k F^k(A^c)$.

\item  For $(x,i), (y,j) \in A$, $(x,i) \prec (y,j)$ if $i = j$ and $x < y$.

\item  $\pi \colon A \rightarrow A$ is defined by $\pi((x,i))= F(x,i)$.
\end{itemize}
Note that if all critical points of $F$ are periodic, then $\pi$ is a bijection. Furthermore the m.c.d. does not depend on the decomposition (up to isomorphisms between m.c.d).
\end{rem}

\begin{defi} Let $f$ be a  multimodal map of type $n$ and consider $F$ an extended map induced by a decomposition $(f_1,\dots,f_n)$ of $F$. If $P$ is a periodic interval for $F$ of period $k$, then we can associate the following m.c.d. $\sigma=<A,\prec,A_c, \pi,c>$
\begin{itemize}

\item $A = \{ F^{i}(P) \colon 0 \leq i < k \}$.

\item $A^c = \{ F^{i}(P) \colon c \in F^{i}(P) \text{ for some critical point $c$ of $F$ }\}$. $c = P$.

\item  For $F^{i}(P) \prec F^{j}(P)$, if $i = j \text{ mod }n $ and $F^{i}(P) < F^{j}(P)$ in the usual order in the real line.

\item  $\pi \colon A \rightarrow A$ is define by $\pi(F^{i}(P))= F^{i + 1 \text{ mod } n}(P)$.
\end{itemize}
Note that $\pi$ is transitive. Furthermore the $\sigma$ does not depend on the decomposition (up to isomorphism between m.c.d).
$\sigma = \sigma(P,f)$ will be called the \textbf{combinatorial type} of the periodic interval $P$. If $P$ is the restrictive interval of the renormalization $Rf$ then $\sigma$ is the combinatorial type of the renormalization of $f$.
\end{defi}

The set of \textbf{$n$-admissible} combinatorial types $\Sigma^n$ is the set of m.c.d. $\sigma$  which are  realized by a $C^1$ renormalizable multimodal map of type $n$. It is easy to see that this is the same that the set of m.c.d. which can be realized by critically finite multimodal map of type $n$ or the set of m.c.d $\sigma$ such that $\# \{ [x] \colon x \in A_\sigma \} = \# A^c_\sigma =n$. Denote by  $\Sigma_k^n$ the subset of $n$-admissible combinatorial types with period bounded by $k$ ($\#A_\sigma \leq kn$).
  
If $F$ is an extended map and $x$ is a point (or an interval) in $I_n$ with a pure itinerary $\omega=\ell_0 \cdots \ell_k$, then $sgn(\omega)$ says if $F^{k+1}$ preserves ($sgn(\omega) = 1$) or reverses ($sgn(\omega) = -1$) the orientation in  $x$. Let $\sigma = <A,\prec,A^c, \pi,c>$ be a transitive m.c.d. and $x \in A$. Let $i \geq  0$ be minimal so that $\pi^{-i}(x) \in \pi(A^c)$. Define $sgn(\sigma)(x)  =sgn(\ell_0(\pi^{-i}(x))\cdots\ell_{i-1}(\pi^{-i}(x)))$. 

\begin{defi} Let $\sigma_1 = <A_1,\prec_1,A^c_1, \pi_1,c_1>$, $\sigma_2 = <A_2,\prec_2,A^c_2,\pi_2,c_2>$ be  m.c.d such that $\pi_1$ is transitive and $\# A^c_1 = \# A^c_2$. The \textbf{product} between $\sigma_1$ and $\sigma_2$ will be  the m.c.d $\sigma = \sigma_2 \ast \sigma_1 = <A,\prec,A^c, \pi,c>$ defined by
\begin{itemize}
\item $A = \{ (x,y)  \colon \Pi_1(x)=(\Pi_1 \circ \pi_1)^ic_1 \text{ and }y \in [\pi_2^ic_2], \text{ for some i}\}$. Moreover $c = (c_1,c_2)$.

\item $(x,y) \prec (\tilde{x},\tilde{y})$ in the following cases:
\begin{itemize}
\item $x \prec_1 \tilde{x}$;
\item  $x = \tilde{x}$, $y \prec_2 \tilde{y}$ and $sgn(\sigma_1)(x)=1$;
\item  $x = \tilde{x}$, $\tilde{y} \prec_2 y$ and $sgn(\sigma_1)(x)=-1$. 
 \end{itemize}

\item $A^c =  \{ (c,\tilde{c}) \in A \colon  c \in A_1^c, \tilde{c} \in A_2^c\}$.

\item $\pi$ is defined by
\begin{equation}
\pi(x,y) = \begin{cases} (\pi_1(x),y)&  \text{ if $x \in A_1 \setminus A^c_1$},  \\
                          (\pi_1(x),\pi_2(y))&  \text{ if $x \in A^c_1$}. \end{cases}
\end{equation}
\end{itemize}
Note that $\# A^c = \# A^c_1 = \# A^c_2$. A m.c.d $\sigma$ is \textbf{primitive} if $\sigma$ does not have a non trivial decomposition $\sigma = \sigma_2 \ast \sigma_1$. 
\end{defi}

\begin{rem} We are primarily interested in to define the $\ast$-product when $\pi_1$ is transitive, but we can give a more general definition when $\# \{(\Pi \circ \pi)^{i}(c_1)\colon i \in \mathbb{N} \} = \# A^c_2$.
\end{rem}

\begin{pro} The $\ast$-product has the following properties:
\begin{itemize}
\item Let $f$ be a $N$ times renormalizable multimodal map of type $n$ such that $f$ is $m$ times renormalizable and $R^m f = A_{P^m} \circ f^{N_m} \circ A_{P^i}^{-1}$. Then the order in the real line of the intervals of disjoint interior in the family $\{f^j(P^m) \colon  j < N_m\}$ are determinated by $\sigma_i \ast \dots \ast \sigma_1$, where $\sigma_i$ is the combinatorial type of the renormalization of $R^{i-1}f$.
\item \label{cfren} Let $\sigma_1$, ..., $\sigma_k$ be an arbitrary sequence of primitive n-admissible m.c.d., where $\sigma_i$ is transitive for $i < k$. Then if $f$ is a critically finite multimodal map of type $n$ with combinatorial type $\sigma = \sigma_k \ast \dots \ast \sigma_1$, then $f$ is $k$ times renormalizable, and the renormalization of $R^{i-1}f$ has type $\sigma_i$.
\end{itemize}
\end{pro}
\begin{proof} The first statement it is easy. Let's to prove the second one. It is sufficient to prove that if $\sigma_1$ and  $\sigma_2$  are $n$-admissible m.c.d., where $\sigma_1$ is transitive, then any critically finite multimodal map of type $n$ with combinatorial type $\sigma_2 \ast \sigma_1$ has a restrictive interval of combinatorial type $\sigma_1$. Let $F$ be an extended map for $f$. Consider the representation of $\sigma$ given by the definition of $\ast$-product. Then $A_\sigma \subset \{(x,y) \colon x \in A_{\sigma_1} \text{ and } y \in A_{\sigma_2}\}$. In other appropriate representation (given in remark \ref{defcf}) $A_\sigma = \{F^{i}c \colon $c$ \text{ is a critical point of $F$} \}$. Let $\phi$ be the unique isomorphism with maps the first representation to the second one. For $x \in A_1$, let $J_x$ be the minimal interval in $I_n$ which contains all points in $\phi(\{ (x,y) \colon (x,y) \in A_\sigma\})$. Then $J_{c_1}$ is a periodic interval for $f$ with combinatorial type $\sigma_1$.
\end{proof}

\subsection{Full families} \label{exist}

\begin{defi} A \textbf{good family} $f_\lambda$ of multimodal maps of type $n$ is a family of multimodal maps of type $n$ such that  $f_\lambda = f_n(\lambda, \cdotp) \circ \dots \circ f_1(\lambda, \cdotp)$, where  the $C^2$-smooth functions $f_i\colon \Lambda \times I \rightarrow I$, and furthermore
\begin{itemize}
\item $f_i(\lambda, \cdotp)\colon I \rightarrow I$ is an unimodal map such that zero is its critical point.

\item The function $(f_1(\cdotp,0), \dots , f_n(\cdotp,0)) \colon \Lambda \rightarrow I^n$ is a homeomorphism.
\end{itemize}
\end{defi}

Let $\sigma = <A,\prec,A^c,c>$. Denote by $\lambda = \lambda(v_1,\dots,v_n)$ be the parameter $\lambda \in \Lambda$ such that $(f_1(\lambda,0),\dots, f_n(\lambda,0))=(v_1,\dots,v_n)$. Let $c_1 = c$, $c_{i+1} \in A^c$ be the unique critical point such that $c_{i+1} \in [\pi(c_i)]$. If $c_i  \in [a]$ define
$$ f^{-1}_{[a],\ell,\lambda} = g_{i, \ell, \lambda} $$
Here $\ell = R$ or $L$ and $g_{i,L, \lambda}$ (resp. $g_{i,R, \lambda}$)  is the orientation preserving (resp. orientation reversing) inverse branch of $f_i(\lambda,\cdotp)$. Consider the compact convex set:
$$ K_\sigma = \{x \in [-1,1]^{A} \colon \text{ if } a_1, a_2 \in A \text{ and } a_1 \prec a_2 \text{ then } x_{a_1} \leq x_{a_2}  \}  $$
Note that $ x \in \partial K_\sigma$ iff $x_{a_1}=x_{a_2}$, with $[a_1]=[a_2]$, $a_1 \neq a_2$. Furthermore, if $x \in K_\sigma$
$$ dist(x,\partial K_\sigma) \leq  C d(x,\partial K_\sigma) \colon =  min_{[a]=[b]} |x_a - x_b|  $$
Here $dist$ is the usual metric defined for a norm in $\mathbb{R}^A$. Define $T \colon int K_\sigma \rightarrow int K_\sigma$ as $T(x)=y$ where $y$ satisfies:
$$ y_a = f^{-1}_{[a],\ell_0(a), \lambda}(x_{\pi (a)})  $$
Here $\lambda = \lambda(x_{\pi(c_1)}, \dots ,x_{\pi(c_n)})$.

\begin{lem} If $ x \rightarrow \partial K_\sigma$ then 
$$ \frac{|T(x) - x|}{dist(x,\partial K_\sigma)} \rightarrow \infty  $$ 
\end{lem}
\begin{proof} The proof  is exactly the proof of lemma 4.1 in pg. 126 of \cite{MS}. We will omit the details. The argument is by contradiction. Assume that
\begin{equation} \label{eqa}  \frac{|T(x) - x|}{d(x,\partial K_\sigma)} \leq K  \end{equation}
Denote $y=T(x)$. Using the argument as in \cite{MS}, we can prove
\begin{enumerate}
\item For $[a]=[b]$:
$$ |x_{\pi (a)} - x_{\pi (b)}| \geq |y_{\pi (a)} - y_{\pi (b)}| - 2K d(x,\partial K_\sigma)$$ 
\item We also have, for $[a]=[b]$:
$$|x_{\pi (a)} - x_{\pi (b)}| \leq C |y_a - y_b| $$   
\end{enumerate}
It follows
$$ |y_{\pi (a)} - y_{\pi (b)}| \leq  C|y_a - y_b| + 2K d(x,\partial K_\sigma)$$
Apply this inequation recursively to obtain 
$$ |y_{\pi^s (a)} - y_{\pi^s (b)}| \leq  C_s |y_a - y_b| + K_s d(x,\partial K_\sigma)$$
Select $a$ and $b$ such that $d(y,\partial K_\sigma) = |y_a - y_b|$ and $s$ such that there is a critical point $\tilde{c} \in A^c$ such that $\pi^s(x_a) \preceq \tilde{c} \preceq \pi^s(x_b)$. Then
$$C_s d(y,\partial K_\sigma) \geq |y_c - y_{\pi (a)}| - K_s d(x,\partial K_\sigma)$$
Because $d(x,\partial K_\sigma) , d(y,\partial K_\sigma) \rightarrow 0$ and since $y_c=0$ is a critical point to the extended map $F_\lambda$, one gets
$$\frac{|y_c - y_{\pi (a)}|}{|x_{\pi(c)} - x_{\pi^2 (a)}|} \rightarrow \infty$$
But
$$\frac{d(y,\partial K_\sigma)}{d(x,\partial K_\sigma)} \geq \frac{1}{C_s} \frac{|y_c - y_{\pi (a)}|}{|x_{\pi(c)} - x_{\pi^2 (a)}|} - \frac{K_s}{C_s} \rightarrow \infty    $$
Which is a contradiction with eq. \ref{eqa}.\end{proof}

\begin{pro} Let $\sigma$ be a  $n$-admissible essential combinatorial type. Then any good family contains  a critically finite multimodal map of type $n$ with combinatorial type $\sigma$.\end{pro}

\begin{proof} By de Melo-van Strien fixed point theorem (see the apendice), there exists a fixed point to the operator $T$ associated to a good family $f_\lambda$ and the essential combinatorial type  $\sigma$.
\end{proof}

\begin{cor}\label{infexists} Let $f_\lambda$ be a good family. Given an infinity sequence of primitive, transitive m.c.d $(\sigma_1,\sigma_2,\dots)$, there exists $\lambda \in \Lambda$ such that $f_\lambda$ is an infinitely renormalizable map with this combinatorial type. 
\end{cor} 

\section{Spaces of polynomials}\label{pol}

\subsection{Polynomials of type n} Consider the polynomials of degree $2^{n}$ such that the dominant coefficient is $1$ and $0$ is a critical point to it. This space can be identified with the $(2^{n}-1)$-dimensional space $P_{n}$ of free coefficients. We say that $p \in P_{n}$ is a \textbf{polynomial of type n} if $p = P_{a_{n}} \circ \dots \circ P_{a_{1}}$, with $P_a(z)=z^2 + a$. Denote this set $Pol_n$.

\begin{pro}\label{unico} $Pol_{n}$ is a complex submanifold of $P_{n}$ with global parameterization 
$$(a_{1},\dots,a_{n}) \rightarrow P_{a_{n}} \circ \dots \circ P_{a_{1}}$$ \end{pro}
\begin{proof}
The following statement, proved by induction, is sufficient to prove the lemma:
Let $\sum_{i \leq 2^{n}} b_{i}x^{i}= P_{a_{n}} \circ \dots \circ P_{a_{1}}$: if $i > 2^{n}-2^{j}$ then $b_{i} = V_{i}(a_{1},\dots,a_{j-1})$. If $i = 2^{n}-2^{j}$ then $b_{i} = C_{j}a_{j} + V_{i}(a_{1},\dots,a_{j-1})$, where $V_i$ are multi-variable polynomials and $C_j \neq 0$.
\end{proof}

The \textbf{connectivity locus} of a family $f_\lambda$, $\lambda \in \Lambda$, of polynomial (or polynomial-like) maps is the set of parameters $\lambda$ such that the filled-in Julia set of $f_\lambda$ is connect. The following result are contained in the stronger results about centered monic polynomials proved by Branner and Hubbard (\cite{BH}). But in our setting the proof is simple:

\begin{pro} The connectivity locus $\mathcal{C}_n$ of $Pol_{n}$ is compact. Moreover all the connected filled-in Julia sets are contained in an uniform neighborhood of zero.  \end{pro}
\begin{proof} We claim that the connectivity locus is contained in the set $\{ P_{a_{n}} \circ \dots \circ P_{a_{1}} \colon |a_{i}|<4 \}$. Indeed, take a polynomial $F=P_{a_{n}} \circ \dots \circ P_{a_{1}}$ in $A_{n}$ outside this set. Let $a_{M}$ such that  $|a_{M}|=max\{|a_{1}|,\dots,|a_{n}| \}$. Consider a critical point $c$ such that $P_{a_{M}} \circ P_{a_{M-1}}\dots \circ P_{a_{1}}(c) = a_{M}$. We claim that $F^{n}(c)$ goes to infinity. This is consequence of a simple fact: if  $b$ is  such that $|b| \geq max\{4,|a_{1}|,\dots,|a_{n}| \}$  then $|b^{2} + a_{i}| \geq 2|b|$ for all $i$. Let $F = P_{a_{n}} \circ \dots \circ P_{a_{1}}$ any polynomial.  Take $b$ like above. Then, using the fact above,  $K(F) \subset B_{|b|}(0)$. In particular, in the connectivity locus the Julia sets are in a fixed neighborhood  of zero. Now is easy to see that the connectivity locus is closed.
\end{proof}

 \begin{pro}\label{muito} For any $N > 1$ there exists $\delta(N)$ with the following property: Consider  $p \in \mathcal{C}_n$ and suppose that there exists $z$ such that $z, p^N(z) \in B_\delta$. Then $B_{2\delta}$ is contained in a periodic component of $K(p)$ which contains a periodic attractor.
\end{pro}
\begin{proof} Let $P = P_{a_n} \circ \dots \circ P_{a_1}$ be a polynomial of type $n$ in $\mathcal{C}$. In particular $|a_i|\leq 4$. Then for any $\delta$ and $N$ there exists $C_1 = C(N,\delta)$ such that $|D(p^{N})(p(z))| \leq C_1$ for all $z \in B_\delta$, $p \in \mathcal{C}$. Furthermore there exists a constant $C_2 = C_2(\delta)$ such that $|p'(z)| \leq C_2 |z|$. Suppose that $z_0, p^N(z_0) \in B_\delta$. For $ z\in B_{2\delta}$, we obtain
$$  |p^N(z)| \leq |p^N(z) - p^N(z_0)| + |p^N(z_0)| \leq  (2 C_1(2\delta,N-1) C_2(2 \delta) \delta + 1)\delta $$
Thus if $\delta$ is small enough, then the map $p^N\colon B_{2 \delta} \rightarrow B_{2 \delta}$ is a strict contraction, and the lemma follows.
\end{proof} 

\begin{cor}\label{juliaab} There exists a constant $C$ so that for any $p \in \mathcal{C}_n$
$$ 1/C \leq diam K(p) \leq C. $$  
\end{cor}

Define the following subsets of $P_{n}\colon$
\begin{itemize}

\item \textbf{$P_{n}^{nd}$} will be the set of $F \in P_{n}$  s.t. $F$ has $2^{n}-1 $ distinct critical points. 'nd' means 'no degenerate'. $P_{n}^{nd}$ is an open set in $P_{n}$.

\item  \textbf{$B_{n}$} will be set of $F \in  P_{n}^{nd}$  s.t $F$ has only $n$ critical values and there exists a  partition $L_{0}, \dots L_{n-1}$ of the critical set of $F$ s.t.  $L_{0} = \{ 0 \}$, $L_{i}$ has $2^{i}$ elements and $F(L_{i})$ is an atomic set. 

\item \textbf{$Pol_{n}^{nd}$}  are the polynomials in $Pol_{n} \cap P_{n}^{nd}$ with  $n$ critical values. Observe that $Pol_{n}^{nd} \subset B_{n}$.
\end{itemize}

Decide if a polynomial of degree $2^n$ belongs to  $Pol_{n}^{nd}$ can be, a priori, difficult:  we would like to take a decision without to look for a explicit decomposition for the polynomial. For $n=2$ we have: 

\begin{pro}  If $n=2$ then $B_{n} = Pol_{n}^{nd}$. \end{pro}

\begin{proof} Take $F \in B_{n}$.  Let $c_1,c_2,0$ be the distinct critical points, and $v_1, v_2$ the critical values, with $F(c_i)=v_2$ and $F(0)=v_1$. Consider the polynomial $Q  = P - v_2$. Observe that $Q(c_i)=0$ and $Q'(c_i)=0$, $i=1,2$. We conclude that $Q = (x -  c_1)^2(x - c_2)^2$. Since that $P'(0) = Q'(0)$ we obtain $0 = -2c_2(c_1)^2 -2(c_2)^2c_1$. This implies $c_1 = - c_2$ (because $c_1$ and $c_2$ are distinct of zero). Hence $F = (x^2 - (c_1)^2)^2 + v_2$.
\end{proof}

It is very easy decide if a polynomial belongs or not to $B_{n}$. For large $n$ the situation is more complex, but sufficiently satisfactory:

\begin{pro}\label{submani} $B_{n}$ is a complex submanifold of $P_{n}$ with dimension  $n$. In particular $Pol_{n}^{nd}$ is an open subset of $B_{n}$.\end{pro}

\begin{proof}
 Each polynomial $F \in P_n$ has $2^{n}-1$ distinct critical points. Order in an arbitrary way the critical points of $F\colon$ $C(F)=\{ c_i  \}$. Define, for P in a neighborhood $V$ of $F$, $c_{i}(P)$, the closest critical point of $P$ to $c_i$. The functions $c_{i}$ are analytic. Define $\phi \colon V \rightarrow {\mathbb{C}}^{2^{n}-1}$ by $\phi(P) = (P(c_{1}(P)),\dots,P(c_{2^{n}-1}(P)))$. We claim that $\phi$ is a local diffeomorphism. Indeed, if $F = \sum_{1 \neq i \leq 2^{n}} a_{i}x^{i}$, then the derivatives of $\phi$ are $\frac{\partial \phi_{j}}{\partial a_{k}}= c_{j}^{k} +\frac{\partial c_{j}}{\partial a_{k}}F'(c_{j})=c_{j}^{k}$ (in this formula, $0^0=1$). Because the critical points are distinct, the Jacobian of $D\phi$ is not zero. Reduce $V$, if necessary, to assume that $\phi$ a diffeomorphism. Take $F \in B_{n}$. Then there exists only one partition of the critical points 
$L_{0}, \dots L_{n-1}$ such that $L_{i}$ has $2^{i}$ elements and $F(L_{i})$ is atomic. Moreover, if $G \in B_{n}$ is close to $F$, $F$ and $G$ respect the same partition. In other words, $F(c_{j})=F(c_{l})$  if and only if $G(c_{j}(G))= G(c_{l}(G)) $. In ${\mathbb{C}}^{2^{n}-1}$, consider the afinne space 
\begin{equation}
S = \{(v_1,\dots,v_{2^{n} -1})\colon \text{$v_i=v_{k} $ iff $F(c_{i})=F(c_{k})$}   \}
\end{equation}
Observe that $S$ is $n$-dimensional and if $G$ is close to $F$,  $G \in B_{n}$ iff $\phi(G) \in S$.  Hence $\phi^{-1}(S\cap \phi(V)$ is a n-dimensional manifold.
\end{proof}

So the previous result will be used as a  'local" criteria to decide if a polynomial belongs to  $Pol_{n}^{nd}$.

\begin{rem} The set $Pol_{n}^{nd}$ is clearly closed in $B_{n}$. Moreover $Pol_{n}^{nd}$ is connect, since it is equal to $Pol_{n}$ up a proper analytic subset. Hence $Pol_{n}^{nd}$ is a connect component of $B_{n}$. It is not difficult to see that $Pol_{n}^{nd}$ is a proper subset of $B_{n}$ for $n > 3$.
\end{rem}

\section{Polynomial like maps}\label{pol-like}

We say that $f\colon U \rightarrow V$, where   $U$ and $V$ are simply connected domains such that  $U$ is compactly contained in $V$, is a \textbf{polynomial like map} if $f$ is a holomorphic ramified covering. The filled-in Julia set $K(f)$ of $f$ is the set of points in $U$ for which all iterates of $f$ are defined. We assume that the McMullen's topology  in the  space of polynomial-like maps (see \cite{Mc1}) is familiar to the reader. Sometimes it is useful work with \textbf{germs} of polynomial like maps: two polynomial like maps $f_i\colon U_i \rightarrow V_i$, $i =1, 2$ define the same germ if
\begin{itemize}
\item The filled-in Julia sets $K(f_1)$ and $K(f_2)$ are equal,
\item the maps $f_1$ and $f_2$ are equal in a neighborhood of $K(f_1)$.
\end{itemize} 

\subsection{Hybrid class} We say that two polynomial like maps $f$ and $g$ are \textbf{hybrid conjugated} if there exists a quasiconformal map $\phi$ defined in a neighborhood of the filled-in Julia set of $f$ and with values in a neighborhood of the filled-in Julia set of $g$ such that $\phi \circ f = g \circ \phi$ and $\overline{\partial} \phi = 0$ in $K(f)$. Now it is a classic result in the  complex dynamics the theorem of Douady and Hubbard (\cite{DH}) which asserts that any polynomial like map is hybrid conjugated with a polynomial. Moreover if $K(f)$ is connect then this polynomial is unique up conjugacies by affine maps. 

The following simple modification of the result of Douady and Hubbard will be a useful tool in the study of the polynomial like contrapart of the concept of multimodal map of type $n$.

\begin{pro}[Straightening lemma] Let $f \colon U_{1} \rightarrow U_{n+1}$ be a polynomial like map, which has the form $f = f_{n} \circ \dots \circ f_{1}$, where $f_{i} \colon U_{i} \rightarrow U_{i+1}$ are ramified coverings of degree $N_i$ and $U_{i}$ are simply connected. Assume that the critical values of $f$ are contained in  $U_{1}$. Then $F$ is hybrid conjugated with a polynomial in the form $P_n \circ \dots \circ P_1$, where $P_i$ is a polynomial of degree $N_i$.\end{pro}

\begin{proof}
First of all, we can assume, using the uniformization Riemann mapping, that $U_{i}= D_{0}(r_{i}) = \{x \colon |x| < r_{i}    \}$, for $i>1$, and  $U_{1} \subset U_{2} \subset U_{3} \dots \subset U_{n+1} $. Assume that the diameter of $U_{n+1}$ is very big. Hence $f_{i} \colon U_{i} \rightarrow U_{i+1}$ are polynomial like maps. We will obtain quasiregular extensions $\tilde{f_{i}}$ of $f_{i}$ and $\tilde{f}$ of $f$ whose are \textbf{compatible}$\colon$  $\tilde{f} =  \tilde{f_{n}} \circ \dots \circ \tilde{f_{1}}$. Choose $\epsilon$ small and define $h_{n+1} = id$ and $A_{n+1}=\tilde{A}_{n+1}=  A(r_{n+1}-\epsilon,r_{n+1})$. Suppose that we had defined $h_{i}\colon A_{i} \rightarrow \tilde{A_{i}}$, where $A_{i}$ is a very fine ring such that the extern boundary of $A_{i}$ is exactly the boundary of $U_{i}$, $\tilde{A_{i}}$ is the pre-image of $A_{n+1}$ by $x^{N_n N_{n-1} \dots N_i}$ and $h_{i}$ is a analytic homeomorphism. Define $A_{i-1}= f^{-1}_{i-1}(A_{i})$ and $h_{i-1}$ as an analytic homeomorphism such that the following diagram commute
\begin{equation}
\begin{CD}
A_{i-1}      @>h_{i-1}>>        {\tilde{A}}_{i-1}\\
@V{f_{i-1}}VV                   @VV{x^{N_{i-1}}}V\\ 
A_{i}        @>>h_{i}>          {\tilde{A}}_{i}
\end{CD}
\end{equation}

Let $H$ be a quasiconformal map which glues the maps $h_{i}$, extending the map to identity outside $U_{n+1}$. Now, we are able to define the quasiregular extensions 
\begin{equation}
  {\tilde{f}}_{i}(x) =\begin{cases}
    f_{i}(x)&      x \in U_{i}, \\
     H^{-1}\circ Q_{N_i} \circ H(x)&       \text{in other case.}
 \end{cases}
\end{equation}
\begin{equation}
  \tilde{f}(x) =\begin{cases}
    f(x)&      x \in U_{1}, \\
     H^{-1} \circ Q_{N_n \dots N_1} \circ H(x)&    \text{in other case.}
 \end{cases}
\end{equation}
Here $Q_{n}(x) \colon = x^{n}$. It is easy to see that these extensions are compatible. Make the pullback of the trivial Beltrami field outside $U_{n+1}$  by the quasiregular mapping $\tilde{f}$. We obtain an invariant Beltrami field $\mu$ for $\tilde{f}$ (defining the field trivial under K(f)). Define $\mu_{i} = ({\tilde{f}}_{n} \circ \dots \circ {\tilde{f}}_{i})_{\ast}\mu$. These Beltrami fields are trivial in a neighborhood of infinity. Let $L_{i}$ be the quasiconformal map so that $\frac{\partial L_{i}}{\overline{\partial} L_{i}} = \mu_{i}$,  $L_{i}(0) = 0$, $L_{i}(\infty) = \infty$, $L_{i}'(\infty) = 1$. Define $L_{n+1}=L_{1}$. Then $P_{i} =  L_{i+1} \circ  {\tilde{f}}_{i}  \circ L_{i}^{-1}$ are polynomials. Moreover $P_n \circ P_{n-1} \dots \circ P_{1} = L_1^{-1} \circ \tilde{f} \circ L_1$ is hybrid conjugated to $f$.
\end{proof}

\begin{figure}
\centering
\psfrag{e}{$h_1$}
\psfrag{f}{$h_2$}
\psfrag{g}{$h_3$}
\psfrag{a}{$f_1$}
\psfrag{b}{$f_2$}
\psfrag{c}{$x^{N_1}$}
\psfrag{d}{$x^{N_2}$}

\includegraphics[width=0.80\textwidth]{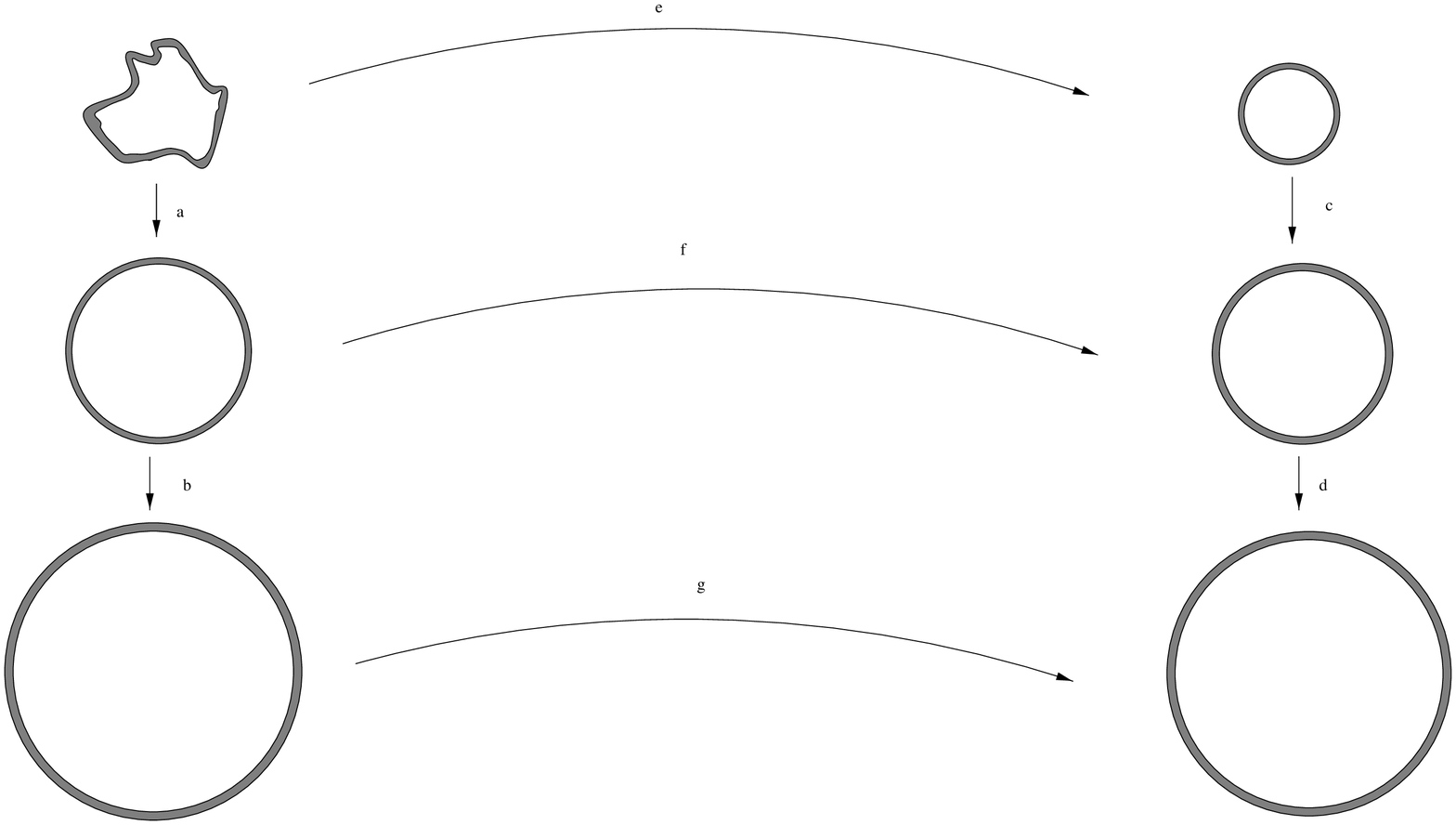}
\caption{}
\end{figure}

 Now we are going to define the polynomial like analogous to the concept of multimodal map of type $n$:

\begin{defi} We say that $f\colon U \rightarrow V$ is a \textbf{polynomial like map of type $\mathbf{n}$} if there exist simply connected domains $U=U_1, \dots, U_n, U_{n+1}= V$ and holomorphic maps $f_i\colon U_i \rightarrow U_{i+1}$, $1 \leq i \leq n$ satisfying
\begin{itemize}  
\item $f_i\colon U_i \rightarrow U_{i+1}$ is a ramified covering maps of degree two,
\item $f = f_n \circ \dots \circ f_1$.
\end{itemize}

By the straightening lemma, any polynomial like map of type $n$ is hybrid conjugated with a polynomial in the form $P_{a_n} \circ \dots \circ P_{a_1}$, where $P_{a}(x) \colon = x^{2}+a$. 
\end{defi}

\begin{rem} Note that we can assume, by the Riemann mapping lemma, that $U_2, \dots, U_{n-1}$ are equal to $\mathbb{D}$.
\end{rem}

We say that a polynomial like map $f\colon U_0 \rightarrow U_n$ of type $n$ is \textbf{real} if there exists a decomposition $(f_1, \dots, f_n)$, $f_i \colon U_i \rightarrow U_{i+1}$, satisfying:
\begin{itemize}
\item The domains $U_i$ are symmetric which respect to the real axis and $I \subset U_i$,

\item The maps $f_i$ preserves the interval $I$. Furthermore $f_i\colon I \rightarrow I$ is an unimodal map such that $(f_1,\dots,f_n)$ is a decomposition for  the multimodal map of type $n$ $f\colon I \rightarrow I$. 
\end{itemize}

\subsection{Extern class} Let $f \colon U \rightarrow V$ be a polynomial like map of degree $d$ and connected filled-in Julia set. Consider $\phi \colon \overline{\co} - K(f) \rightarrow \overline{\co} - \overline{\mathbb{D}}$ be the Riemann mapping such that $\phi(\infty)=\infty$. Then the map $ g = \phi \circ f \circ \phi_{-1}\colon \phi(U - K(f)) \rightarrow \phi(V - K(f))$ is defined in an open set $A-\overline{\mathbb{D}}$, where $A$ is a neighborhood of $\partial \mathbb{D}$. We can invert $g$ along $\partial \mathbb{D}$ to obtain a holomorphic map $g\colon \tilde{U} \rightarrow \tilde{V}$ defined in a neighborhood of $S^1$. Here $\tilde{U}$  is the union of $\overline{\phi(U - K(f))}$ with its inversion along $S^1$. It is easy to see that $g$ is an expanding map of degree $d$. The map $g$ is called the \textbf{extern map} of $f$. Note that $g$ is defined up to affine conjugacies. Indeed, the external map can be defined when $K(f)$ is not connected, but this will not be used here.

\begin{pro}\label{commutes} Let $A$ and $B$ be neighborhoods of $S^1$ and  $h \colon A - \overline{\mathbb{D}} \rightarrow B -\overline{\mathbb{D}}$ be a homeomorphism  which commutes with $x^d$. Then $h$ has a continuous extension to $A - \mathbb{D}$ such that $h(z)=\alpha z$ in $S^1$, with $\alpha^{d-1}=1$.\end{pro}
\begin{proof} Let 
$$ A_R = \{z \in \co \colon 1 < |z| < R  \}$$
and   
$$ B_R = \{z \in \co \colon z = x + y\cdot i, \text{ with } x, y \in \mathbb{R}, 0 < |z| < ln R    \}$$
Denote $Q(z)=z^d$ and define   $\tilde{Q} \colon B_R \rightarrow B_{ R^d}$ by $\tilde{Q}(z)= d \cdot z$ and  $\mu(z) = e^{2\pi z}$. The diagram
\begin{equation}
\begin{CD} 
B_R     @>\tilde{Q}>>       B_{R^d}\\
@V\mu VV                   @VV\mu V\\ 
A_R        @>>Q>          A_{ R^d}
\end{CD}
\end{equation}
commutes. Fix an arbitrary $R > 1$ and assume, without loss of generality, that the domain of $h$ is $A_{r}$, for some $r < R$ and $B -\overline{\mathbb{D}}$ is contained in $A_R$. Consider the fundamental annulus $A = Q^{-2}(A_{R^d} - A_R)$ and a compact set $\tilde{A} \subset B_R$ such that $\mu(\tilde{A})=A$. Then for any point $\tilde{z} \in \mu^{-1}(Q^{-2}(A_{R^d}))$, there exist $j, i \in \mathbb{N}$ such that $   \tilde{T}^{i}(\tilde{z}) + j \in \tilde{A}$.

Select a homeomorphism $\tilde{h}$ such that $\mu \circ \tilde{h} = h \circ \mu$. The transformation $\tilde{h}$ is defined in the open set $\tilde{A} \cap A_R$ and with values in $\tilde{B} \cap A_R$, where $\tilde{A}=  \mu^{-1} (A - \overline{\mathbb{D}})$ and $\tilde{B}=  \mu^{-1} (B - \overline{\mathbb{D}})$. Since $h \circ T = T \circ h$, one gets
\begin{equation}\label{equac}
 \tilde{h} \circ \tilde{T} = \tilde{T} \circ \tilde{h} + k
\end{equation}
For some $k \in \mathbb{N}$ and for points in $A_{r^{1/d}}$. Since 
\begin{equation}
 (\tilde{h}+j) \circ \tilde{T} = \tilde{T} \circ (\tilde{h}+j) + k + j \cdot (d-1)
\end{equation}
we can assume, replacing $\tilde{h}$ by an appropriate translation of $\tilde{h}$ by an integer, that $0 \leq  k < d -1$. Apply the equation \ref{equac} recursively to obtain
$$ \tilde{h} \circ \tilde{T}^{-i} = \tilde{T}^{-i} \circ \tilde{h} - k \cdot  d^{-i} \cdot \frac{d^i -1}{d -1 },$$
for $j > 0$ and points in $A_{r^{1/d}}$. Let 
\begin{equation}\label{sobreD} D = sup_{z \in \tilde{A}} dist_{B_R}(h(z), z - \frac{k}{d-1})
\end{equation}
Then for any $ \tilde{x} \in \mu^{-1}(Q^{-2}(A_{R^d})$, select $i, j$ such that  $x = \tilde{T}^{-i}(z) - j$, with $z \in \tilde{A}$. If $\tilde{x} \in \tilde{A}$ then 
$$dist_{B_R}(\tilde{h}(\tilde{x}), \tilde{x} - \frac{k}{d-1} )$$ 
$$ = dist_{B_R}(\tilde{h}(\tilde{T}^{-i}(z) - j), \tilde{T}^{-i}(z) - j - \frac{k}{d-1})$$
$$ =  dist_{B_R}(d^{-i} \cdot \tilde{h}(z) - k \cdot d^{-i} \cdot \frac{d^{i} -1}{d -1 } - j, d^{-i} \cdot z - j - \frac{k}{d-1} )   $$
$$ =  dist_{B_R}(d^{-i} \cdot (\tilde{h}(z) - k \cdot \frac{d^{i} -1}{d -1 } - j \cdot d^{i}), d^{-i} \cdot (  z - k \cdot \frac{d^i}{d-1} -  j \cdot d^{i}))   $$
$$ \leq   dist_{B_R}(\tilde{h}(z) +  \frac{k}{d -1 } - k \cdot \frac{d^{i}}{d -1 } - j \cdot d^{i},  z  - k \cdot \frac{d^i}{d-1}  -   j \cdot d^{i})   $$
$$ =  dist_{B_R}(\tilde{h}(z), z -  \frac{k}{d -1 }) \leq D,  $$
The above proof is a variation of the Douady and Hubbard's proof(\cite{DH}) when $k = 0$. So $dist_{A_R}(h(x),\alpha \cdot x) \leq D$, for $\mu(\tilde{x})=x$, because $\mu \colon B_R \rightarrow A_R$ is a local isometry. Hence  $\alpha \cdot x - h(x) \rightarrow 0$  in the Euclidean topology, when $x \rightarrow S^1$, since $\rho_{A_R} \geq C \frac{1}{dist(x,S^1)}$ near to $S^1$.
\end{proof}

We will denote $[h]=\alpha$. Notice $[h_1 \circ h_2]= [h_1] \cdot [h_2]$, for two homeomorphisms $h_i$ which commutes to $x^d$.

\begin{defi} Let $f\colon U \rightarrow V$ and $g\colon \tilde{U} \rightarrow \tilde{V}$ be polynomial-like maps. Let $h_{1}\colon U-K(f) \rightarrow \tilde{U}-K(g)$ and $h_{2}\colon U-K(f) \rightarrow \tilde{U}-K(g)$ be two homeomorphisms such that $ g =  h_{i} \circ f \circ h_{i}^{-1}$. Then $h_{1}\circ h_{2}^{-1}$ is an automorphism of $f$ in $U-K(f)$. Choose a homeomorphisms  $\phi\colon U-K(f) \rightarrow A(1,r)$ such that $\phi \circ f = T \circ \phi$, with $T(x)=x^d$ and  $d = deg \ f$. Then $\phi \circ h_{1}\circ h_{2}^{-1}  \circ \phi^{-1}$ is an automorphism of $T$. Define 
$$[f,g;h_1,h_2]=[\phi \circ \psi  \circ \phi^{-1}]$$
Observe that $[f,g;h_1,h_2]$ is well defined, since if $\phi_1$ and $\phi_2$ are two conjugacies between $f$ and $T$ then
$$ [\phi_1 \circ \psi  \circ \phi_1^{-1}] = [\phi_1 \circ \phi_2^{-1}]\cdot [\phi_2 \circ \psi  \circ \phi_2^{-1}] \cdot [(\phi_1 \circ \phi_2^{-1})^{-1}] = [\phi_2 \circ \psi  \circ \phi_2^{-1}] $$ 
\end{defi}
The number $[f,g;h_1,h_2]$ was introduced by Douady and Hubbard(\cite{DH}) to study when it is possible to glue conjugacies:

\begin{cor}[Gluing conjugacies:\cite{DH}] Let $f$ and $g$ be polynomial like maps, let $h_1 \colon V_f  \rightarrow V_g $ and $h_2 \colon V_f - K(f) \rightarrow V_g - K(g)$ be conjugacies. If $[f,g;h_1,h_2]=1$ then there exists a conjugacy $h \colon  V_f  \rightarrow V_g$ such that 
\begin{itemize}
\item The map $h$ coincides with  $h_1$ in $K(f)$.

\item The map $h$ coincides with $h_2$ in $V_f - K(f)$.
\end{itemize} 
\end{cor}

\begin{cor} Let $f\colon U \rightarrow V$ and $g\colon \tilde{U} \rightarrow \tilde{V}$ be polynomial like maps with the same hybrid and extern class. If there is an extern equivalence $h_{1}$ and a hybrid equivalence $h_{2}$ such that $[f,g,h_{1},h_{2}]=1$, then $f$ and $g$ are affine conjugated.
\end{cor}

\begin{pro}\label{homotopy} Let $H \colon [0,1] \times A -\overline{\mathbb{D}}\rightarrow B -\overline{\mathbb{D}}$ be a map such that $A$ and $B$ are neighborhoods of $S^{1}$ and  $h_t = H(t,\cdot)$ are homeomorphisms which commutes with $x^d$. Then there is a continuous extension $H \colon [0,1] \times A -\mathbb{D}\rightarrow B -\mathbb{D}$. In particular $h_t(z) = \alpha \cdot z$, for all $t$, where $\alpha^{d-1}=1$. 
\end{pro}
\begin{proof} Using notations as in the proof of proposition \ref{commutes}, we have that, since the maps $h_t$ are homotopic, that $ \tilde{h}_t \circ \tilde{T} = \tilde{T} \circ \tilde{h}_t + k$, where $k$ does not depend on $t$. Moreover the constant $D$ in equation \ref{sobreD} can be select independent of $t$, which is sufficient to prove the lemma.
\end{proof}

\begin{pro} Let $f$ and $g$ be polynomial like maps with connected Julia sets  and let $h\colon U_f \rightarrow U_g$ be a topological conjugacy between $f$ and $g$. Then there is an isotopy $H \colon [0,1] \times \tilde{U}_f \rightarrow \co$ such that:
\begin{itemize}
\item For each $t$, $h_t=H(t,\cdot)$ is a conjugacy between $f$ and $g$ which coincides with $h$ in $K(f)$,
\item $h_0 = h$,
\item $h_1$ is quasiconformal outside $K(f)$.
\end{itemize} 
\end{pro}
\begin{proof} As in the proof of lemma 3.2 in \cite{PR99}, we can construct a  
isotopy $h_t \colon V_f - K(f) \rightarrow  V_g - K(g)$ so that $h_0 = h$, $h_t$ is a conjugation and furthermore $h_1$ is quasiconformal. Consider $L(t,x)=h_0^{-1}\circ h_t(x)$. For each $t$, $L(t,\cdot)$ is an automorphism of $f$ in $V_f - K(f)$. It is well know that there is a quasiconformal map $\phi$ which is a conjugacy between $f$ in $V_f - K(f)$ and $x^d$, where $d = deg \ f$. So $\tilde{L}(t,x) = \phi (L(t,\phi^{-1}(x)))$ satisfies the hypothesis of the proposition \ref{homotopy}, and furthermore $\tilde{L}(0,x)=x$. In particular, $\tilde{L}(t,x)$ is at a finite hyperbolic distance of $x$, independent of $t \in [0,1]$. Since $\phi$ is quasiconformal,  the same is true for $L(t,x)$ with respect the hyperbolic metric of $V_f - K(f)$. In particular we can extend $L(t,x)$ in a continuous way to $K(f)$ setting $L(t,x)=x$, for $x \in K(f)$. Extend $H$ setting $H(t,x) = h_0(L(t,x))$. 
\end{proof}

The existence of a conjugacy which is quasiconformal outside $K(f)$ and coincides to $h$ in $K(f)$ was carry out in a more general case by Przytycki and Rohde(\cite{PR99}),  and for hyperbolic maps by McMullen and Sullivan (\cite{
SM98}). For quadratic-like maps, Lyubich (\cite{Lyu}) proves the existence of an isotopy (in the quadratic case, the isotopy is the identity in $K(f)$).

\subsection{Compact sets}\label{compact}

In a neighborhood of the locus  of connectivity of polynomials of type $n$, select a holomorphic moving fundamental annulus $A_p$. This means that for each $p \in \mathcal{C}_n$ there exist a  neighborhood $\Lambda$ of $p$ and a map $\psi \colon \Lambda \times A_p \rightarrow \co$ such that
\begin{itemize}
\item For each $\tilde{p} \in \Lambda$, $\psi(\tilde{p},A_p)=A_{\tilde{p}}$.
\item For each point $z$ in the annulus $A_p$, $\psi(\cdot,z)$ is a holomorphic function.
\item If $U_p$ and $V_p$ are respectly bounded simple connected domains whose boundaries are the internal and external boundaries of the annulus $\psi(p,A_p)$, then $p\colon U_p \rightarrow V_p$ is a polynomial-like map. 
\end{itemize} 

To find such $\phi$, let $C$ be a large circle centered in zero which contains all the Julia sets for polynomials of type $n$ contained in the locus of connectivity. Then, for $p$ near to $\mathcal{C}$, $p^{-1}(C)$ is a Jordan curve contained in the disc whose boundary is $C$. We obtain a polynomial-like representation for $p$. Furthermore the set $C \cup p^{-1}(C)$ moves holomorphicaly in a neighborhood of $p$. Let $\psi$ be this holomorphic motion. This holomorphic motion can be extended for all points in $\co$. Then we define $A_{\tilde{p}}$ as the annulus delimited by  $C$ and  $\tilde{p}^{-1}(C)$. We will fix this holomorphic moving annulus on the polynomials of type $n$ in the rest of this paper. 

\begin{defi}\label{de} We say $f\colon U \rightarrow V$,  a polynomial like map of type $n$ with connect Julia set, belongs to $\mathcal{P}_n(C)$ if  there are  a $C$-quasiconformal map $\phi\colon \co \rightarrow \co$ and $p \in \mathcal{C}_n$  such that 
\begin{itemize}
\item $\phi(V_p - U_p)= V - U$,
\item $   \phi \circ p = f \circ \phi$ in $U_p$.
\end{itemize}
\end{defi}

\begin{pro} The set $\mathcal{P}_n(C)$ is compact up to affine conjugacies.
\end{pro} 
\begin{proof} Let $f_i$ be a sequence in  $\mathcal{P}_n(C)$.  Replacing $f_i$ by a polynomial like map which is affine conjugated to it, we can assume $diam K(f_i) = 1$. Consider  $C$-quasiconformal maps  $\phi_i$ and polynomial maps $p_i$ as in definition \ref{de}. Since $\mathcal{C}_n$ is compact, select a subsequence, if necessary, such that $p_i \rightarrow p \in \mathcal{C}_n$. Since the Julia set of $p \in \mathcal{C}_n$ has the diameter away of infinity and zero, and $\phi_i(K(f_i))=K(p_i)$, selecting a  subsequence we can assume that  $\phi_i$ converges to a $C$-quasiconformal map. It is not difficult to see that $\phi \circ p \circ \phi^{-1}$ is an analytic map in $\phi(V_p - U_p)$.  
\end{proof}

\begin{cor}\label{delta} For any $\delta > 0$ there exists $n = n(C,\delta)$ so that if $f\colon U \rightarrow V$ belongs to $\mathcal{P}_n(C)$ then $f^{-n}(V) \subset \delta \text{-}K(f)$.
\end{cor}
\begin{proof} Easy. \end{proof}

\begin{lem}[\cite{Mc2}]\label{comp1} Let $f\colon U \rightarrow V$ be a polynomial like map  with connect filled-in Julia set and $\epsilon$-$K(f) \subset V$. Then the germ of $f$ has a representation $f\colon \tilde{U} \rightarrow \tilde{V}$ such that:
\begin{itemize}
\item The boundaries of $\tilde{U}$ and $\tilde{V}$ are $C(\epsilon)$-quasicircles,
\item $diam \tilde{V} \leq \tilde{C}(\epsilon) diam K(f)$,
\item $mod(\tilde{V} - \tilde{U}) \geq m(\epsilon)$.
\end{itemize}
\end{lem}

When $p$ is a polynomial of type $n$ with connect Julia set, we can select a polynomial-like restriction $p\colon \tilde{U}_0 \rightarrow \tilde{U}_n$ with the above properties in the following way: Let $\phi\colon \co - K(f) \rightarrow \co - \mathbb{D}$ be the Riemann map such that $\phi(\infty)=\infty$. We have 
$$\phi \circ p (x) = (\phi(x))^{d}$$
for $x \in \co - K(p)$ and $d = deg  p$. Let $\mathbb{D}_r = \{z \in \co \colon |z| \geq r  \}$. Define
$$U_0 = \phi^{-1}(\mathbb{D}_{\exp \frac{m}{d-1}} - \overline{\mathbb{D}}) \cup K(p) \text{ and } U_n = \phi^{-1}(\mathbb{D}_{\exp \frac{dm}{d-1}} - \overline{\mathbb{D}}) \cup K(p)$$
Then $mod (U_{n+1} - U_0) = m$. It is easy to prove that $\partial \tilde{U}_0$ and $\partial \tilde{U}_n$ are $C(m)$-quasicircles. To prove that $diam \tilde{U_n} \leq \tilde{C}(m) diam K(p)$, recall that the diameter of $K(p)$ is bounded above and below, by lemma \ref{juliaab}. So it is sufficient to prove that $diam U_{n+1} \leq C(m)$. Indeed, consider the Green function $G(x)= log |\phi(x)|$. Since $\mathcal{C}$ is compact, for any $\epsilon > 0$ there exists $R_\epsilon$ such that for $|z| \geq R_\epsilon$ and $p \in \mathcal{C}$ one have
$$  \frac{1}{d} \log (1 - \epsilon) + \log  |z| \leq G(x) $$ 
( e.g., see the proof of Theorem 2.1 in \cite{FS}), which clearly implies $diam U_{n+1} \leq C(m)$, since $G(\partial U_{n+1})= dm/(d-1)$. The advantage of this polynomial like restriction is that the annulus $U_{n+1} - K(f)$ and $U_0 - K(f)$ are invariant by the extern automorphisms 
$$\phi^{-1} \circ  R_\alpha   \circ \phi \colon \co - K(f) \rightarrow \co - K(f)$$ 
where $R_\alpha(x)= \alpha x$ and $\alpha^{d-1}=1$.

Denote by $\mathcal{P}_n(C,\tilde{C},m)$ the set of polynomials like maps of type $n$ which admits a decomposition $f = f_n \circ \dots \circ f_1 \colon U_1 \rightarrow U_{n+1}$, $f_i\colon U_i \rightarrow U_{i+1}$, such that
\begin{itemize}
\item The filled-in Julia set $K(f)$ is connect;
\item The boundaries of $U_i$ are either $C$-quasicircles, for $i=1,n+1$, or the unit disc, otherwise; 
\item $diam U_{n+1} \leq \tilde{C} diam K(f)$;
\item $mod(U_{n+1} - U_1) \geq m$;
\item The critical point of $f_i$ is 0.
\end{itemize} 

\begin{pro}\label{comp2} Let $f\colon U_1 \rightarrow U_n$ and $\tilde{f}\colon  \tilde{U}_1 \rightarrow \tilde{U}_n$ be polynomial like maps in of type $n$ which belongs $\mathcal{P}_n(C,\tilde{C},m)$ which are conjugated by a one-to-one continuous map  $h_0$ in a neighborhood of their filled-in Julia sets. Then there exists a one-to-one continuous map  $h_1\colon \co \rightarrow \co $ between $f$ and $\tilde{f}$ with the following properties:
\begin{itemize}
\item The map $h_1$ is a conjugacy: $\tilde{f}  \circ h_1 = h_1 \circ f$ in $U_1$
\item $h_1(U_n - U_1)=\tilde{U}_n - \tilde{U}_1$.
\item $h_1 = h_0$ in $K(f)$.
\item $h_1$ is $C(C,\tilde{C},m)$-quasiconformal in $\co - K(f)$.
\end{itemize}
In particular, if $h_0$ is a hybrid conjugacy then $h_1$ is a hybrid conjugacy.
\end{pro}
\begin{proof} Assume that $diam K(f) = diam K(\tilde{f})=1$. Note that, since the boundary of $U_i$ is a quasicircle, the map $f_i$ has a quasiregular extension in a neighborhood of $\overline{U_i}$. Indeed, let $\alpha_i\colon \co \rightarrow \co$ be a $C(C_1)$-quasiconformal map which is conformal in $\mathbb{D}$ and  maps $\mathbb{D}$ in $U_i$. Then $\alpha_{i+1}^{-1} \circ f_i \circ \alpha_i \colon \mathbb{D} \rightarrow \mathbb{D}$ extends to a rational map $g_i$ which is a expansive map of degree $2$ in $S^1$. Hence $\alpha_{i+1} \circ g_i \circ \alpha_i^{-1} \colon \co \rightarrow \co$ is a regular map of degree two in a neighborhood of $\overline{U}_i$. The same can be done for $\tilde{f}_i$. Let $\phi_n\colon \co \rightarrow \co$ be a $C_n(C)$-quasiconformal map which maps  $U_n$ in $\tilde{U}_n$. Since $f_n(0)$ and $\tilde{f}_n(0)$ are contained in the Julia sets of $f$ and $\tilde{f}$, these points are at a definitive Euclidean distance of $\partial U_{n+1}$ and $\partial \tilde{U}_{n+1}$. Thus, by lemma \ref{quasi}, we can assume that  $\phi_n(f_n(0))=\tilde{f}_n(0)$.  By induction, suppose that we have constructed   $C_j(C,m)$-quasiconformal maps $\phi_i\colon U_i \rightarrow \tilde{U}_i$, for $j$ between $i$ and $n$, such that
\begin{itemize} 
\item $\phi_{i+1} \circ f_{i} =  \tilde{f}_{i} \circ \phi_{i}^{-1}$ in $\partial U_{i}$.
\item $\phi_{i}(f_{i}(0))=\tilde{f}_i(0)$.
\item The quasiconformality of $\phi_i$ is bounded by a constant $C_i(C,m)$.
\end{itemize} 
Let  $\phi_{i-1}\colon U_{i-1} \rightarrow \tilde{U}_{i-1}$ be a lift of  $\phi_{i}$ (in other words: $\phi_{i} \circ f_{i-1} =  \tilde{f}_{i-1} \circ \phi_{i-1}^{-1}$), which has the same quasiconformality than   $\phi_{i}$. Because the critical values of $f$ does not intercept $U_{n+1} - U_{1}$, and the modulus of this annulus is bounded below, $f_{i-1}(0)$ is at a bounded distance of $0$ in the hyperbolic metric on $U_{i-1}$. Since the same can be said about $\tilde{f}_{i-1}(0)$,  by lemma \ref{quasi}, if necessary  modify  $\phi_{i-1}$ such that $\phi_{i-1}(f_{i-1}(0)) = \tilde{f}_{n-1}(0)$, and additionally the new $\phi_{i-1}$ is $C_{i-1}(C,m)$-quasiconformal. In particular $\phi_{n} \circ f = \tilde{f} \circ \phi_0$ in $\partial U_0$. We can find a $C(C,m,M)$-quasiconformal map $H \colon \co \rightarrow \co$ such that (1) $H$ is equal to $\phi_{n+1}$ outside $U_{n+1}$ (2) $H$ is equal to $\phi_1$ in $\overline{U_1}$. Hence $H$ is a quasiconformal map which maps $\overline{U_{n+1} - U_1}$ in $\overline{\tilde{U}_{n+1} - \tilde{U}_1}$ and conjugates $f$ and $\tilde{f}$ in the boundary of this fundamental annulus. Now, with the usual pullback argument, construct a $C(C,m,M)$-quasiconformal conjugacy $H$ between $f$ and $\tilde{f}$ in $U_{n+1} - K(f)$ in $\tilde{U}_{n+1} - K(\tilde{f})$ such that $f(U_{n+1} - U_1)= \tilde{U}_{n+1} - \tilde{U}_1$. For the last step, to obtain a conjugacy which extends to $K(f)$, the result follows of the particular case when $\tilde{f}$ is a polynomial and the annulus $\tilde{U}_{n+1} - \tilde{U}_1$ is invariant by the extern automorphisms of $\tilde{f}$. Select a extern automorphism $R\colon \co - K(\tilde{f}) \rightarrow \co - K(\tilde{f})$ so that $[h, R \circ H; f, \tilde{f}]=1$. Thus the conjugacy $h_1$ in $K(f)$ glue with the external conjugacy $R \circ H$ and the new map $h_1$ has the same quasiconformality of $H$ outside $K(f)$.   
\end{proof}

\section{Renormalization}\label{ren} 
\subsection{Infinitely renormalizable polynomials} Here we will work with a more natural parameterization of polynomial of type n. We will consider the family  of polynomials in the form $f= f_{a_1} \circ \dots \circ f_{a_n}$, where $f_a(x) = -2a x^{2} + 2a - 1$, $a \in \mathbb{\co}$. Note that $\alpha f(1/ \alpha)$ belongs to $Pol_n$, with
$$ \alpha^{2^n -1 }= - 2^{2^n -1} a_1 a_2^2 a_3^4 \dots a_n^{2^{n-1}}. $$ 

For each $f$ in this family there is at most $2^n -1$ polynomials in $Pol_n$  affine conjugated to it. Furthermore, if $a_i \in \mathbb{R}$ then there is exactly an real map in $Pol_n$ affine conjugated to $f$. Because the results of section \ref{exist}, there is, for each combinatorial type $\sigma=(\sigma_1,\sigma_2, \dots)$, at least one infinitely renormalizable multimodal map of type $n$ in this family   which has type $\sigma$. Denote the set of infinitely renormalizable real polynomials of type $n$ by $Pol_n^\infty$ and the subset of  $Pol_n^\infty$ with $C$-bounded combinatoric by $Pol_n^\infty(C)$. Denote by $F$ the extended map associated to a decomposition in quadratic polynomials  $(f_{a_1},\dots, f_{a_n})$.  The bounded geometry of the postcritical set of $f \in Pol_n^{\infty}(C)$ is a particular case of a result proved by J. Hu\cite{HU98} (see also J. Hu's  thesis\cite{H95}). The following result can be proved as in \cite{Sm} (for notation, see the end of the introduction):
  
\begin{pro}[Bounded geometry]\label{bounded} Let $q,r,s$ be arbitrary critical points of $F$ so that $Q_{-i}^{k+1}$ and $R_{-j}^{k+1}$ are contained in $S^k_{-\ell}$. The following quantities are $C_1(C)$-commensurable:
\begin{itemize}
\item The lengths of $Q_{-i}^{k+1}$, $R_{-j}^{k+1}$ and $S^{k}_{-\ell}$,
\item The distance between  $Q_{-i}^{k+1}$ and $\partial S^{k}_{-\ell}$,
\item The distance between $Q_{-i}^{k+1}$ and $R_{-j}^{k+1}$, if these intervals do not touch.  
\end{itemize}
\end{pro}

Denote by $\mathcal{P}_n^\infty(C_1,C_2)$ the set of maps in $\mathcal{P}_n(C_2)$ which are hybrid conjugated with polynomials in $Pol_n^\infty(C_1)$. The main technical result in renormalization theory is 

\begin{pro}[Complex bounds:\cite{Sm}]\label{complex} Let $f$ be a map in $Pol_n^\infty(C_1)$.  Then there exist $k_0(C_1)$ and $C_2(C_1)$ so that any renormalization $R^k(f)$, $k \geq k_0$, has a polynomial-like extension $R^k(f)\colon U \rightarrow V$ in $\mathcal{P}^{\infty}_n(C_1,C_2)$. Furthermore, the renormalization is unbranched: $P(f) \cap V = P(f) \cap K(R^k(f))$.    
\end{pro}
\begin{proof} Here we have a very nice situation: the map f is a polynomial, it has negative Schwartzian derivative and moreover satisfies properties analogous to the standard conditions (see \cite{Sm}). It is easy  verify in the proof of the complex bounds for analytic multimodal maps \cite{Sm} that $k_0$ and $C_2$(for this use lemmas \ref{comp1} and \ref{comp2}) can be select independent of $f \in Pol_n^{\infty}(C_1)$.
\end{proof}

J. Hu(\cite{H95}) stated a complex bounds result for bimodal maps in the Epstein class and bounded combinatorics, but the outline of the proof seems to be incomplete.

\begin{pro}\label{cont} The following statements holds:
\begin{itemize}
\item Let $f_n$ be maps in $Pol_n^{\infty}(C)$ with combinatorics $\sigma_n$. If $\sigma_n$ converges to $\sigma$ then any limit $f_\infty$ of a  subsequence of $f_n$ has combinatorics $\sigma$.
\item The postcritical set moves continuously in $Pol_n^\infty(C)$.
\end{itemize}
\end{pro}
\begin{proof} Let $P^k_n$ be the restrictive interval associate to the $k$-th renormalization of $f_n$, $0 \in P^k_n$. Since the period of $P^k_n$ is bounded by $p_0(k,C)$, the length of $P^k_n$ can not be small, otherwise $f_n$ will contain a periodic point  which attracts zero, which is absurd. So we can assume that $P^k_n$ converges to a periodic interval $P^k_\infty$ for $f_\infty$, which proves that $f_\infty$  belongs to  $Pol_n^\infty(C)$. In particular,  all periodic points of $f_\infty$ in $I$ are reppeling (because non reppeling periodic points attracts a critical point), so the periodic point in the boundary of $P^k_n$ converges to a periodic point in $\partial P^k_\infty$. Thus $P^k_n$ is the unique restrictive interval associate to the $k$-th renormalization of $f_\infty$ and the $k$-th restrictive interval moves continuously in $Pol_n^\infty(C)$ and so do the postcritical set.    
\end{proof}

\subsection{Renormalization for  polynomial-like maps}Let $f\colon U \rightarrow V$ be a polynomial like map. A \textbf{pre renormalization} of $f$ is a polynomial like map  $g\colon \tilde{U} \rightarrow \tilde{V}$ such that
\begin{enumerate}
\item $\tilde{U} \subset U$,
\item  $g = f^i$ for some $i > 0$,
\item The filled-in Julia set $K(g)$ is connect.
\end{enumerate}
Note that a pre-renormalization of a polynomial-like map of type $n$ is a polynomial like map of type $k$, for some $k \geq 0$. This is a consequence of the following observation: if $g_1$ and $g_2$ are holomorphic maps such that $g_1 \circ g_2 \colon U \rightarrow V$ is a proper map, then $g_1 \colon U \rightarrow g_1(U)$ and $g_2 \colon g_1(U)\rightarrow g_1(U)$ are  proper maps.

\begin{lem}\label{juliacom} Let $f \colon U \rightarrow V$ be a polynomial like map in $\mathcal{P}_n(C)$ with a pre renormalization $g=f^m\colon \tilde{U} \rightarrow \tilde{V}$ such that $c \in K(g)$, where $c$ is the critical point mapped to zero by conjugacies with a polynomial of type $n$. Then $diam K(g) \geq C_1(C,m) diam K(f)$ 
\end{lem}
\begin{proof} Follows of lemma \ref{muito}. 
\end{proof}

\begin{lem}\label{po} Let $g_1=f^n\colon U_1 \rightarrow V_1$ and $g_2=f^n\colon U_2 \rightarrow V_2$ be pre renormalizations of a polynomial like map $f$. Consider $K = K(g_1) \cap K(g_2)$. Then one of the following statements holds:
\begin{enumerate}
\item $K = \phi$.
\item $K = \{p\}$, where $p$ is a repelling periodic point of $f$.
\item  $K$ is the filled-in Julia set of a pre-renormalization $g\colon \tilde{U} \rightarrow \tilde{V}$ of $f$. Moreover  $deg(g) \leq min\{deg(g_1), deg(g_2)\}$ and the equality holds iff $K(g)=K(g_1)$ or $K(g)=K(g_2)$. 
\end{enumerate}
\end{lem}  
\begin{proof} Follows of the connectedness principle by McMullen (pg. 90 in \cite{Mc1}) that $K = K(g_1) \cap K(g_2)$ is connect. Let $\tilde{U}$ be the connect component  of $U_1 \cap U_2$ which contains $K$. Then $g=f^n \colon \tilde{U} \rightarrow f^n(\tilde{U})$ is a polynomial like map, and moreover $K(g)=K$, since $K$ is totally invariant by $g$. Hence we obtain item 2, if $deg \ g=1$, or 3, otherwise. The last statement of item 3 follows of lemma 5.11 in \cite{Mc1}.
\end{proof}

\begin{rem}\label{re} A special case is when  each critical point of the extended map $F$ associate to $f$ can be accumulated by points in the closure of the postcritical orbit of $c$. In this case if $g_i=f^k\colon U_i \rightarrow V_i$, $i = 1, 2$ are two pre renormalizations with $\deg g_i \leq 2^n$ then or $K(g_1)=K(g_2)$ and $\deg g_i = 2^n$  either $K(g_1)\cap K(g_2)$ is at most a repelling period point. In particular if $g$ is a pre-renormalization of $f$ whose domain contains  $c$  and it has degree at most $2^n$ then $\deg g = 2^n$ and any pre renormalization of $g=f^k\colon U \rightarrow V$, for fixed $k$, whose domain contains $c$ define the same germ of polynomial like map of type $n$. Then we can call this germ $g$ as the \textbf{renormalization} of $f$.   
\end{rem}

We do not know if there is a canonical way of define renormalization when $f$ do not satisfy the hypothesis in the previous remark. However, in the case of real polynomial like map  $f\colon U \rightarrow V$ of type $n$, we can to use external rays which arrive in the boundary points of the restrictive interval $P$ to find a degenerate polynomial like extension to the renormalization. After modify the domain near to the boundary points of P, we obtain:

\begin{pro} Let $f\colon U \rightarrow V$ be  a real polynomial like map of type $n$  which is renormalizable in in the sense of section \ref{multimodal}. Let  $P$ be the restrictive interval associate with the renormalization $Rf$ and let $k$ be minimal such that $f^{k}(P) \subset P$. Then there exists a pre renormalization  $\tilde{g}\colon \tilde{U} \rightarrow \tilde{V}$ of degree $2^{n}$ such that $K(\tilde{g}) \cap \mathbb{R} = P$.  
\end{pro}

\begin{rem} If $g_1$ and $g_2$ are two pre renormalizations of degree $2^n$ such that the restrictive interval $P$ is contained in $K(g_i)$. Then $K(g_1)\cap K(g_2)$ is the filled in Julia set of a pre-renormalization $g$. Since $K(g_1)\cap K(g_2)$ contains $P$, $deg(g)= 2^n$. By lemma \ref{po} $K(g)=K(g_1)=K(g_2)$. Thus  $g_i$ define the same polynomial like germ. This germ will be called the \textbf{complex renormalization} of $f$. 
We say that a polynomial like map of type $n$ is renormalizable if it is hybrid conjugated with a renormalizable real polynomial like map of type $n$. Denote by  $\mathcal{P}_n^{\infty}(C_1,C_2)$  the subset of maps in $\mathcal{P}_n(C_2)$ which are hybrid conjugated with a real infinitely renormalizable polynomial of type $n$ with combinatorics bounded by $C_1$. 
\end{rem}

Let $K$ be a closed set in $\co$. We say that $K$ has \textbf{$C$-bounded geometry} if $1/C < \sup_{A \in \mathcal{A}} \mod A < C$ where $\mathcal{A}$ is the set of annulus $A \subset \co - K$ such that the both components of $\co - A$ contain points in $K$.   
The following result will be used a lot of times:

\begin{pro}\label{propinf} For $f \in \mathcal{P}^\infty_n(C_1,C_2)$, the followings holds:
\begin{enumerate}
\item $K(f)$ has empty interior.
\item For almost every point $x$ in the Julia set,  $f^{n}(x) \rightarrow P(f)$.
\item The small Julia sets touch at most in an unique point.
\item  The postcritical set $P(f)$ and $K(f)$ moves continuously in $\mathcal{P}^\infty_n(C_1,C_2)$.
\item There exists a constant $C(C_1,C_2)$ such that the postcritical set $P(f)$ has $C$-bounded geometry.
\item The $k$-th renormalization, for $k \geq k_0(C)$, has a polynomial-like extension $R^{k}(f)\colon U^k \rightarrow V^k$ which belongs to $\mathcal{P}^\infty_n(C,C_2)$. Here $C=C(C_1,C_2)$. 
\item There exists $j(C_1,C_2)$ so that 
$$ \frac{diam \ K(R^{k+j}(f))}{diam \ K(R^{k}(f))} \leq \frac{1}{2}$$
\end{enumerate} 
\end{pro}
\begin{proof} The proof of 1 is exactly as in the unimodal case (see \cite{Mc1}): we can assume that $f$ is a polynomial. Suppose, by contradiction that $K(f)$ has interior. Then, by the Sullivan's classification of periodic components, the interior of $K(f)$ contains an attractor or a Siegel disc: the first case is impossible because $P(f)$ is a Cantor set (and any attractor attracts a critical point) and the second one does not hold because the boundary of a Siegel disc must be contained in the postcritical set. The second statement is consequence of the ergodic or attract  theorem(\cite{Mc1}). Item 3 is consequence of remark \ref{re}.  Item 4 and 5 follow of the same statements for polynomials (propositions \ref{cont} and \ref{complex}). The last item is obvious for polynomials, since  $diam \ K(R^{k}(f))$ is commensurable with the length of $P_0^k$, which goes exponently fast to zero (proposition \ref{bounded}). Now the general case is easy.\end{proof}

If $M$ is a hyperbolic Riemann surface, denote by $\norm{\cdot}_M$ the hyperbolic metric in $TM$ and $dist_M(\cdot,\cdot)$ the hyperbolic distance. We will denote by $dist(\cdot,\cdot)$ the Euclidian distance.

\begin{pro}\label{distcont} There exists a constant $D$  such that, for any $f\colon U \rightarrow V \in \mathcal{P}_n^\infty(C_1,C_2)$, 
$$dist_{V - P(f)}(z,f^{-1}(P(f))) \leq D$$ 
for $z \in \overline{f^{-1}(V-U)}$.
\end{pro}
\begin{proof} It is easy too see there is a bound for $diam_{V - P(f)}\overline{f^{-1}(V-U)}$ which depends only on $C_1$ and $C_2$. So it suffice to proof that there is a point $z \in \overline{f^{-1}(V-U)}$ whose hyperbolic distance to $f^{-1}(P(f))$ is under control. Firstly, assume that $p\colon U \rightarrow V$ is a polynomial in $Pol_n^{\infty}(C_1)$ and $V_p - \overline{U_p}$ is the holomorphic moving fundamental annulus $A_p$ selected in section \ref{compact}. We will prove that there exists $\tilde{D}$ such that for each $p \in Pol_n^\infty(C_1)$ there exist points $x_p \in V_p - \overline{U_p}$, $y_p \in p^{-1}(P(p))$ and a topological disc $B_p$ such that $x_p, y_p \in B_p$ and $ \overline{B_p} \subset V_p - P(p)$ satisfying:
$$ dist_{B_p}(x_p,y_p) \leq \tilde{D}  $$
Indeed,   for each  $p$ with combinatorics bounded by $C_2$, select a point $z_0 \in \overline{f^{-1}(V_p-U_p)}$ and  a topological disc $B$ which contains $z_0$ and a point  $z_1 \in p^{-2}\{0\} - P(p)$. Furthermore $\overline{B} \subset V_p - P(p)$. since $P(p)$ and $\partial V_p$  moves continuously with $p$, for $\tilde{p}$ close to $p$ we have  $ \overline{B} \subset  V_{\tilde{p}} - P(\tilde{p})$. Furthermore there is a point $\tilde{z}_1$ of $\tilde{p}^{-2}({0}) - P(\tilde{p})$ close to $z_1$. In particular, $dist_B(z,\tilde{z}_1)$ is under control. Since $Pol_n^{\infty}(C_1)$ is compact, the prove is finished.\end{proof}

Let $f\colon U_1 \rightarrow U_{n+1}$ be a polynomial like map of type $n$. For each  decomposition $f_i\colon U_i \rightarrow U_{i+1}$, we can associate the extended map $F\colon \mathcal{U} \rightarrow \mathcal{V}$, where $\mathcal{U} = \{(x,i)\colon x \in U_i, 1\leq i \leq n \}$ and $\mathcal{V} = \{(x,i)\colon x \in U_i, 2\leq i \leq n+1 \}$, defined by
$$F(x,i) = (f_i(x), i+1 \ mod \ n)$$
Thus $F$ is a ramified covering map between the Riemann surfaces $\mathcal{U}$ and $\mathcal{V}$. If $c_i$ is the critical point of $f_i$, define the postcritical set of $F$ by
$$P(F) = \overline{\cup_i \cup_j F^{j}(c_i,i)}$$
Assume now that $f \in \mathcal{P}_n^\infty(C_1,C_2)$. Then $P(F)$ is a Cantor set with bounded geometry. In particular $\mathcal{V} - P(F)$ is a $M(C_1,C_2)$-uniform domain. Let $R^jf=f^{m(j)} \colon U \rightarrow V$ be a renormalization of $f$. We say that $K(R^jf), F(K(R^jf)), \dots, F^{m(j) - 1}(K(R^jf))$ are \textbf{small Julia sets}. Clearly, by complex bounds, each small Julia set is the Julia set of a polynomial like map $g$ in $\mathcal{P}_n^\infty(C_3(C_1,C_2),C_2)$. Furthermore the fundamental annulus of $g$ does not intersect the postcritical set of $F$. For each small Julia set $K$ we can associate the closed geodesic $\gamma$ in the hyperbolic domain $\mathcal{V} - P(F)$ which separes $P(F) \cap K$ and $P(F) - K$. These geodesics cut the domain $\mathcal{V} - P(F)$ in subsets which we will call \textbf{pieces}. 

\begin{lem}\label{geod} Let $F\colon \mathcal{U} \rightarrow \mathcal{V}$ be the extended map defined above. Let $K$ be a small Julia set for $F$, $P = P(F) \cap K$ and let $\gamma$ be the closed geodesic in $\mathcal{V} - P(F)$ which separates $P$ and $P(F) - K$. There exists $C_3$, which depends only on $ C_1$ and $C_2$, so that
\begin{itemize}
\item The hyperbolic diameter of $\gamma$ in $\mathcal{V} - P(F)$ is $C_3$-commensurable to one,
\item $dist(\gamma,P)$ and  $diam \ P$ are $C_3$-commensurable,
\item The Euclidean diameters of $K$, $\gamma$, $P$ and the Euclidean length of $\gamma$  are $C_3$-commensurable.
\end{itemize}  
\end{lem}
\begin{proof} Firstly note that, by lemma \ref{juliacom} and a priori bounds, $K$ is contained in a larger small Julia set $K'$ so that the diameters of $K$ and $K'$ are commensurable. Furthermore the diameters of $K$ and $P$ are also commensurable. The hyperbolic length of $\gamma$ is not large because there is fundamental ring for each renormalization with definitive modulus. Since $\mathcal{V} - P(F)$ is a $M(C_1,C_2)$-uniform domain, by lemma \ref{ere} the Euclidean diameter of $\gamma$ is commensurable to $diam \ P$. It is easy to see that the Euclidean length of $\gamma$ is also commensurable to $diam \ P$. If the hyperbolic length is small then the Euclidean diameter of $P$ will be small relative to $K'$, which is a contradiction. The second statement is consequence of the first one and proposition \ref{ere}.\end{proof} 

We say that a map belongs to $\mathcal{F}(f)$ if the graph of $g$ is contained in $\{(x,y)\colon f^i(x)=f^j(y)\}$, for some $i,j \geq 0$. We are going to prove that there are copies of the small Julia sets close to any point in $J(f)$, in any scale:

\begin{pro}[Small Julia sets everywhere]\label{julia} Let $f \in P^{\infty}_n(C_1,C_2)$. There exist $C_3(C_1,C_2)$ and $C_4(C_1,C_2)$ with the following property: For any $z$ in $J(f)$ and $\alpha \in (0,1)$ there exists a polynomial like map  $g\colon U \rightarrow V$, $g \in \mathcal{F}(f)$ so that
\begin{itemize} 
\item $g \in \mathcal{P}_n^{\infty}(C_1,C_4)$,
\item The diameter of $K(g)$ is $C_3$-commensurable to $\alpha \cdot diam(K(f))$,
\item $dist(z,K(g)) \leq C \cdot \alpha$. 
\end{itemize}
\end{pro}
\begin{proof} The prove is quite similar to the proof in the unimodal case (see \cite{Mc2}): Consider a decomposition in ramified coverings of degree two $f = f_1 \circ \dots \circ f_n$, $f_i \colon U_i \rightarrow U_{i+1}$, and the associate extended map  $F$ defined in $\mathcal{U} = \{(z,i) \colon z \in U_{i} \text{ and } 1 \leq i \leq n \}$ by
$$ F(z,i) = (f_i(z), i + 1 \ mod \ n )   $$
Define the postcritical set of $F$ by $P(F)=\cup_{i} F^i(P(f))$. If $\norm{\cdot}_{\mathcal{V} - P(F)}$ denote the hyperbolic metric on $\mathcal{V} - P(F)$ (extended to $\infty$ on $P(F)$) then 
$$ \norm{F'(x)\cdot v}_{\mathcal{V} - P(F)} > \norm{v}_{\mathcal{V} - P(F)}$$
, with $v$ in the tangent space of $x$, since $F^{-1}(P(F)) \supset P(F)$. Let $v_k = DF^{k}(x)\cdot v$, with $x \in J(F)$. Then 
$$  \norm{v_k}_{\mathcal{V} - P(F)} < \norm{v_{k+1}}_{\mathcal{V} - P(F)} \rightarrow \infty $$ 
, since $\cup_n F^{-n}(P(F))$ is dense in $J(F)$. For $x \in J(F)$, select a vector $v$ in its tangent space so that $ |v| = \alpha$. There are 3 cases:
\begin{itemize}
\item $dist(x, P(F)) \leq \alpha$,
\item $dist(x,P(F)) > \alpha$ and there exists $k$ such that $ \norm{v_k}_{\mathcal{V}  - P(F)} \leq \epsilon$ and  $\norm{v_{k+1}}_{\mathcal{V}  - P(F)} \geq  1/ \epsilon$,
\item $dist(x,P(F)) > \alpha$  and $\norm{v_k}_{\mathcal{V}  - P(F)} \sim 1$.
\end{itemize}
Here $\epsilon$ is sufficiently small so that the McMullen's argument(\cite{Mc2}) works in the second case. The first and second  cases are more easy and we will omit the proof. For details, see \cite{Mc2}. Assume the last situation.  In particular $\norm{v}_{\mathcal{V} - P(F)} < 1$, because $\mathcal{V} - P(F)$ is a $M(C_1,C_2)$-uniform domain. 

Consider the piece, defined by closed geodesics, which contains $F^{k}(x)$. Let $\gamma_j$ be the exterior geodesic and let $\gamma_{1,j+1}, \dots, \gamma_{i,j+1}$ be the interior boundary geodesics. Denote by $P(i,j+1)$ the subset of the postcritical set bounded by $\gamma_{i,j+1}$. Select $\ell$ minimal so that we can do the univalent pullback of the domain $V$ bounded by $\gamma_j$ along the inverse orbit $F^{k}(x), F^{k-1}(x), \dots, F^{\ell}(x)$. This means there exists   a simply connected domain $V'$ satisfying
\begin{itemize}
\item $F^{\ell}(x)\in V'$.
\item The map $F^{k -\ell}$ is univalent in $V'$ and moreover $F^{k -\ell}(V')=V$.
\item The domain $V'$ contains a critical value $v$ of $F$.
\end{itemize}
Denote by $V_i$ the domain bounded by $\gamma_{j+1,i}$ and let $V_i'$ be the corresponding domain in $V'$. Let $\tilde{V} = F^{-1}(V')$ Then 
$$g=F^{k-\ell+1}\colon \tilde{V} \rightarrow V$$
is a proper map of degree two. Since the postcritical set is contained in $\cup_i V_i$, the critical value in $V'$ is contained in some $V_{i_0}'$, for some $i_0$. Choose an arbitrary $V_{i_1}$, $i_1 \neq i_0$.Let $\beta_1$ and $\beta_2$ be two paths inside the piece which contains $F^{k}(x)$ so that:
\begin{itemize}
\item The initial  point of both is $F^{k}(x)$.
\item The end point of both is a point in $\gamma_{j,i_1}$. 
\item The Jordan curve defined by $\beta_1$ and $\beta_2$ is not homotopic to a constant curve in $V - F^{k-\ell}(v)$. 
\item The hyperbolic diameter of $\beta_i$ on $\mathcal{V}  - P(F)$ is bounded. 
\end{itemize}

Let $\tilde{\beta}_1 \cup \tilde{U}^1$ and $\tilde{\beta}_2 \cup \tilde{V}^{2}$ be lifts with respect to $g$ of the simply connected sets $\beta_1 \cup V_{i_1}$ and $\beta_2 \cup V_{i_1}$ so that $\tilde{\beta}_i$ is an arc whose initial point is $F^{\ell-1}(x)$. Note that $\tilde{V}^{1}$ and $\tilde{V}^{2}$ are disjoint and one of them, say $\tilde{V}^{1}$, does not intersect the postcritical set of $F$. So all  inverse branches of $F$ are well defined on  $\tilde{\beta}_1 \cup \tilde{V}^1$. So let $h$ the inverse branch of $F^k$, defined in $\beta_1 \cup V_{i_1}$ so that $h(F^k(x))=x$. Since $\norm{v_k}_{\mathcal{V}  - P(F)} \sim 1$ and $\beta_1 \cup \gamma_{j+1,i_1}$ has bounded hyperbolic diameter in $\mathcal{V}  - P(F)$, we obtain, by corollary \ref{disto},
$$|D h(z)| \sim \frac{|v_0|}{|v_k|}$$ 
for all $z \in \beta_1 \cup \gamma_{j+1,i_1}$. By the maximum principle the same distortion control holds in $U_{i_1}$. There exists a small Julia set $K$ inside $V_{i_1}$ whose diameter is commensurable to $diam \gamma_{j+1,i_1} \sim |v_k|$. By $1/4$-Koebe lemma (use that $dist(\gamma_{j+1,i_1} ,P(j+1,i_1))  \geq C \ diam P(j+1,i_1)$) and the above distortion control, the set $h(K)$ has diameter commensurable with $|v_0|=\alpha$. Moreover $dist(x,h(K)) \leq C \alpha$, which proves the proposition. 
\end{proof}

\begin{figure}
\centering
\psfrag{U}{$U$}
\psfrag{U0}{$U_{i_0}$}
\psfrag{U1}{$U_{i_1}$}
\psfrag{b1}{$\beta_1$}
\psfrag{b2}{$\beta_2$}
\psfrag{tb1}{$\tilde{\beta}_1$}
\psfrag{tb2}{$\tilde{\beta}_2$}

\psfrag{gx}{$F^{k}(x)$}
\psfrag{g0}{$\gamma_{j+1,i_0}$}
\psfrag{g1}{$\gamma_{j+1,i_1}$}
\psfrag{Ut0}{$\tilde{U}^1$}
\psfrag{Ut1}{$\tilde{U}^2$}
\psfrag{Fk}{$F^{k -a +1}$}
\psfrag{x}{$F^a(x)$}
\includegraphics[width=\textwidth]{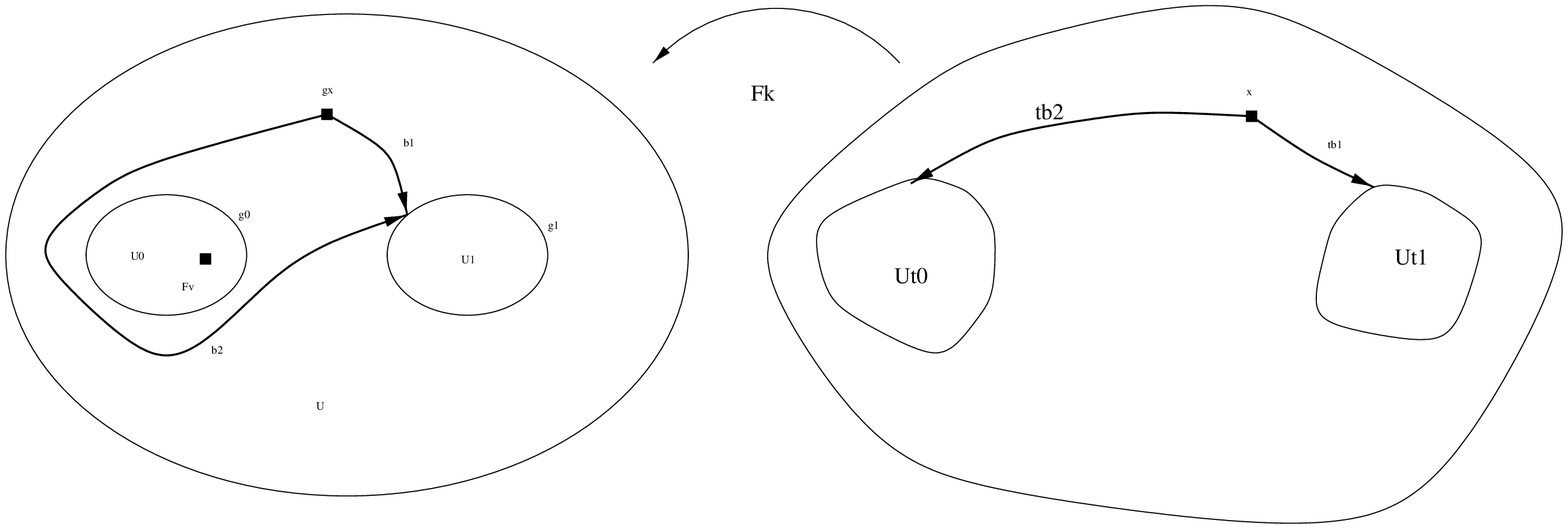}
\end{figure}

\begin{cor}\label{noinvariant} The following holds:
\begin{itemize}
\item A polynomial like map in $\mathcal{P}_n^\infty(C_1)$ does not support invariant line fields in it Julia set,
\item The hybrid class on $\mathcal{P}_n(C_2)$ is continuous at  points in $\mathcal{P}_n^\infty(C_1,C_2)$,
\item The set $\mathcal{P}_n^\infty(C_1,C_2)$ is compact.
\end{itemize}
\end{cor}
\begin{proof}  Suppose, by contradiction, there exists an invariant line field $\mu$ supported on the Julia set $K(f)$ on $f \in \mathcal{P}_n^\infty(C_1)$. Select an almost continuity point $x \in J(f)$ to $\mu$. There are polynomial like maps, with definitive modulus, in all scales around $x$, which preserves $\mu$. After an affine conjugation, we can assume that a subsequence of these polynomial like maps converge to a polynomial like map which preserves a straight line field, which is a contradiction. The second statement is consequence of the lemma in pg.313 of \cite{DH}. The last statement is an immediate consequence of the first ones.
\end{proof}

\section{Hybrid conjugacy}\label{hyb}

As in the unimodal case (see \cite{S} and \cite{MS}) and bimodal case(\cite{H95}), real maps with same bounded combinatorics are hybrid conjugated. To be more precise:

\begin{thm}\label{pullback} Two real polynomial-like maps of type $n$, $f$ and $\tilde{f}$, infinitely renormalizable map with same bounded combinatorics are hybrid conjugated.
\end{thm}

Let $f$ and $\tilde{f}$ be multimodal maps of type $n$ with decompositions $(f_1,\dots,f_n)$ and $(\tilde{f}_1,\dots,\tilde{f}_n)$. If $c_i$ (resp. $\tilde{c}_i$) is  the critical point of $f_i$ (resp. $\tilde{f}_i$), define $v_i = f_n \circ \dots \circ f_i (c_i)$ (resp. $\tilde{v}_i = \tilde{f}_n \circ \dots \circ \tilde{f}_i (\tilde{c}_i)$.

\begin{lem}[Lifts exist]\label{lift} Let $f\colon U_0 \rightarrow U_n$ and $\tilde{f}\colon \tilde{U}_0 \rightarrow \tilde{U}_n$ be real polynomial-like maps of type $n$ with  the same inner itinerary and such that $v_i < v_j$ iff $\tilde{v}_i < \tilde{v}_j$. Then the following holds:  For any continuous bijection  $H_n\colon U_n \rightarrow \tilde{U}_0$, real in the real line and increasing in $\mathbb{R}$, such that $H(f(C(f)))= \tilde{f}(C(\tilde{f}))$, there exists a continuous bijection $H_0\colon U_0  \rightarrow \tilde{U}_0$, real in the real line and increasing in $\mathbb{R}$, such that $H_n \circ f = \tilde{f} \circ \tilde{H_0}$.
\end{lem}
\begin{proof} We will define, by induction, homeomorphisms $H_i \colon U_i 
\rightarrow \tilde{U}_i$, real in the real line and increasing, such that $H_{i+1} \circ \tilde{f}_i = f_i \circ H_{i}$. Assume that we have defined $H_i$. We claim that $H_i(f_{i-1}(c_{i-1}))= \tilde{f}_{i-1}(\tilde{c}_{i-1})$. Indeed, consider $A_i = (f_n \circ \dots \circ f_i)^{-1}(v_{i-1})$ and $\tilde{A}_i = (\tilde{f}_n \circ \dots \circ \tilde{f}_i)^{-1}(\tilde{v}_{i-1})$. Since $H_n \circ f_n \circ \dots \circ f_i = \tilde{f}_n \circ \dots \circ \tilde{f}_i \circ H_i$ and $H_n(v_{i-1})=\tilde{v}_{i-1}$, we have $H_i(A_i)= \tilde{A}_i$. Since $H_i$ is increasing, it suffice to show that if $f_{i-1}(c_{i-1})$ is the $j$-th point in $A_i$, with respect to the order in the real line, then $\tilde{f}_{i-1}(\tilde{c}_{i-1})$ is also the $j$-th point in $\tilde{A}_i$. But this follows of lemma \ref{pos}.3, since $f_{i-1}(c_{i-1})$ and $\tilde{f}_{i-1}(\tilde{c}_{i-1})$ have the same inner itinerary.  

Now we can find a homeomorphism $H_{i-1}\colon U_{i-1} \rightarrow \tilde{U}_{i-1}$, real in the real line and increasing such that 
\begin{equation}
\begin{CD}
U_{i-1}-\{c_{i-1}\}      @>H_{i-1}>>       {\tilde{U}}_{i-1}-\{\tilde{c}_{i-1}\}\\
@V{f_{i-1}}VV                   @VV\tilde{f}_{i-1}V\\ 
U_{i}-\{f_{i-1}(c_{i-1}) \}        @>>H_{i}>          {\tilde{U}}_{i} - \{\tilde{f}_{i-1}(\tilde{c}_{i-1}) \}
\end{CD}
\end{equation}
commutes, which proves the lemma.
\end{proof}

\begin{lem}\label{posconj} If $f$ and $\tilde{f}$ are as in Theorem \ref{pullback}, then there exists a quasiconformal map $h\colon \co \rightarrow \co$, which is  real in the real line and increasing such that $h(P(f))=P(\tilde(f))$ and $h \circ f = \tilde{f} \circ h$ in $P(f)$.
\end{lem}
\begin{proof} This lemma is consequence of the bounded geometry (see proposition \ref{bounded}). For details see \cite{S} or \cite{MS}(last chapter)).
\end{proof}

\begin{proof}[Proof of Theorem \ref{pullback}] Replacing $U, \tilde{U}, V$ and $\tilde{V}$ by smaller domains, we can assume that the boundary of these domains are quasicircles. Using similar arguments as in lemma \ref{comp2}, we can construct 
a quasiconformal map $h_1\colon \overline{V}- U \rightarrow \overline{\tilde{V}}- \tilde{U}$  which conjugates $f$ and $\tilde{f}$ in $\partial U$. Since the maps $f$ and $\tilde{f}$ are real, we can assume that $h_1$ is symmetric with respect to the real line. Let $h_2$ be a symmetric quasiconformal map which conjugates $f$ and $\tilde{f}$ in the postcritical set. Construct a $C$-quasiconformal map $H_0 \colon V \rightarrow \tilde{V}$, for some $C$, which is symmetric,  increasing in the real line, equal to $h_1$ in $\overline{V}- U$ and equal to $h_2$ in a neighborhood of $I$. As $f$ and $\tilde{f}$ have the same combinatorial type, the relative positions of $v_i$ and $\tilde{v}_i$ are the same. Thus we can use lemma \ref{lift}. Furthermore $f$ and $\tilde{f}$ has the same inner itinerary. Define inductly $H_j \colon V \rightarrow \tilde{V}$, a $C$-quasiconformal map symmetric, increasing and  such that $H_j \circ f = f \circ H_{j+1}$. Note that $H_j$ is $C$-quasiconformal conjugacy in the postcritical set and $V - f^{-(j-1)}(U)$. Moreover $H_j = H_{j+1}$ in $V - f^{-(j-1)}(U)$. Since $K(f)$ has empty interior, the sequence $H_j$ has an unique limit $H$, which is a conjugacy between $f$ and $\tilde{f}$. Indeed $H$ is a hybrid conjugacy by lemma \ref{noinvariant}.  
\end{proof}

\begin{thm} The   set $Pol_n^{\infty}(C)$ is a Cantor set.
\end{thm}
\begin{proof} Let $\Sigma$ be the set of primitive, transitive m.c.d with combinatorics bounded by $C$. By lemma \ref{infexists}, for any infinity sequence $(\sigma_1,\sigma_2,\dots)$, with $\sigma_i \in \Sigma$, there exists an infinitely renormalizable real polynomial map of type $n$ with this combinatorial type. By the previous theorem,  two real polynomial maps of type $n$ with same combinatorics are hybrid conjugated and so affine conjugated, since they are polynomials. The point $0$ must be a fixed point for this affine map. Thus the conjugacy must be the identity. Let $\Pi \colon \Sigma^\mathbb{N} \rightarrow A_n$ be the application which maps each sequence   $\alpha = (\sigma_1,\sigma_2,\dots)$ in the unique real polynomial map  $p_\alpha$ of type $n$ with this combinatorial type. If $\alpha_i \in \Sigma^\mathbb{N}$ is a sequence which converges to $\alpha$, then any accumulation point $p$ of the sequence $p_{\alpha_i}$ is a real infinitely renormalizable polynomial  of type $n$ with combinatorics $\alpha$. So $p = p_\alpha$. Hence $\Pi$ is a homeomorphism between the Cantor set $\Sigma^\mathbb{N}$ and $Pol_n^{\infty}(C)$.   
\end{proof}

\section{Convergence of renormalization}\label{tow}

\subsection{Towers}

\begin{defi} A \textbf{ bi-infinity tower} $\mathbf{f}$$=<f_i>_{i \in \mathbb{Z}}$ of type $n$ with parameters $C_1$, $C_2$ and $k$ is a family of polynomial-like maps $f_{i}\colon U_{i} \rightarrow V_{i}$ of type $n$, $i \in \mathbb{Z}$, such that 
\begin{itemize}

\item The maps $f_i$ belongs to $\mathcal{P}_n^{\infty}(C_1,C_2)$;

\item For any $i \in \mathbb{Z}$, there exists $ a \leq k$ so that  $R^a(f_{i-1}) = f_i$.
\end{itemize}  

if furthermore we assume that 

\begin{itemize}
\item If $j \geq i$ then  $V_j$ is contained in $V_{i}$.
\end{itemize}

we say that $\mathbf{f}$ is a \textbf{fine} tower. Denote by $\mathcal{T}_n(C_1,k,C_2)$ (resp. $\mathcal{T}_n^{fine}(C_1, k, C_2)$) the set of bi-infinity towers (resp. fine towers) with parameters $C_1$, $k$,  $C_2$. 
\end{defi}

McMullen(\cite{Mc2})  supply the set of towers with the following sequential convergence: We say that the the sequence of towers $\mathbf{f}_n$ converges to tower $\mathbf{f}_\infty$ if 
\begin{itemize}
\item For any $i \in \mathbb{Z}$ there exists $a$ so that $f_{i+1,n} =  R^{a}(f_{i,n})$, for large $n$;

\item $f_{i,n}$ converges to $f_{i,\infty}$.
\end{itemize}

\begin{pro} The sets $\mathcal{T}_n(C_1, k, C_2)$ and $\mathcal{T}_n^{fine}(C_1, k, C_2)$ are compact.\end{pro}

Select an arbitrary $j_1 \in  \mathbb{Z}$. Then  $A_j = V_j - U_j$, $j_0 < j < j_1$, are disjoint essential annulus in $V_{j_0} - K(f_{j_1})$. Because $diam \ K(f_{j_1}) > 0$ and $\mod A_j > m(C_1)$, $\co = \cup_j V_j$.
 
We say that a line field $\mu$ is invariant by the tower $\mathbf{f}$ if $\mu$ is invariant for each $f_i$, $i \in \mathbb{Z}$.

\begin{pro}[Construing bi-infinite towers]\label{cbit} Let $f_{j,i}\colon U_{j,i} \rightarrow V_{j,i}$; with   $i \in \mathbb{N}$ and $|j| \leq j(i)$, $j(i) \rightarrow _i \infty$; be polynomial-like maps in $\mathcal{P}^{\infty}_n(C_1,C_2)$ such that there exists $k$ satisfying 
$$f_{j+1,i} = R^a(f_{j,i}), \text{ where } a \leq k.$$
Then we can select a subsequence $i_k$ such that $f_{j,i_k} \rightarrow f_{j,\infty}$, where $\mathbf{f}_{\infty}=<f_{j,\infty}>$ is a tower in $\mathcal{T}_n(C_1, k, C_2)$. If $V_{j+1,i} \subset V_{j,i}$, then  $\mathbf{f}_\infty \in \mathcal{T}_n^{fine}(C_1, k, C_2)$.
\end{pro}

\begin{pro}[Construing conjugacies]\label{cconj} Let 
$$f_{j,i} \colon U_{j,i} \rightarrow V_{j,i} \text{ and } \tilde{f}_{j,i} \colon \tilde{U}_{j,i} \rightarrow \tilde{V}_{j,i}$$
with $|j| \leq j(i)$, be as in the previous lemma. Let $h_i \colon \co \rightarrow \co$ be  $k$-quasiconformal maps  so that:
\begin{itemize}
\item $h_i$ is a hybrid conjugacy between $f_{j,i}$ and $g_{j,i}$, for $|j |\leq j(i)$,

\item $h_i(U_{j,i})=\tilde{U}_{j,i}$.
\end{itemize}
Then we can select a subsequence $i_k$ such that 
$$<f_{j,i_k}>_{|j| \leq j(i_k)} \text{ and }  <\tilde{f}_{j,i_k}>_{|j| \leq j(i_k)}$$
converge to bi-infinite towers $\mathbf{f}_\infty$ and $\mathbf{g}_\infty$ and $h_{i_k}$ converges to a conjugacy between these towers.
 \end{pro}

Fix $f = <f_i>_{i \in \mathbb{Z}}$ a bi-infinite tower.

\begin{lem}\label{contra} Suppose that $z \in U_i$ with  $f_i(z) \in V_i - P_i$. Then 
$$\hiper{f_i'(z)}{\co - P}{\co -P} \geq 1$$. 
\end{lem}

\begin{proof} Since $f_i \colon U_i - Q_{i} \rightarrow V_i - P_i$ is a covering map 
$$ \hip{f_i'(z)}{i} = 1,  $$
and furthermore
$$\hip{\incl{i}}{i} < 1,$$
we obtain 
$$\hiper{f'(z)}{V_i  - P_i}{V_i - P_i} > 1.$$
Since the hyperbolic metric in $V_i - P_i$ converges to the hyperbolic metric in $\co - P$, we obtain the lemma.
\end{proof}

\begin{lem}[Strict contraction(\cite{Mc2})]\label{escontra} There exists $\lambda > 1$ with the following property: Let $z \in U_i$ be such that $f_i(z) \in V_i - U_i$. Then 
$$ \norma{f_i'(z)}{\co - P} \geq \lambda  $$  
\end{lem}
\begin{proof} We sketch the McMullen's proof: Consider $ j < i$. We have $f_i = f^{a}_j$, for some $a > 0$. so 
$$\hiper{f_i'(z)}{f^{-a}_j(V_j - P_j)}{V_j - P_j} = 1.$$
Since $z \in f_i^{-1}(V_i - U_i)$, by the inclusion contraction lemma (proposition 4.9 in \cite{Mc2}): 

$$\norm{i_{f^{-a}_j(V_j - P_j), V_j - P_j}}_{U_i - Q_i, V_j - P_j} \leq C(D) <  1$$
It follows that  $\norm{f_i'(z)}_{V_j - P_j} \geq \lambda(D)$. Now it is suffice to observe  that  $\rho_{V_j - P_j} \rightarrow \rho_{\co - P}$.     
\end{proof}

Let $i < j$ be such that $z_0 \in U_j - J(f_i)$. Consider $t_\ell$, with $ i \leq \ell \leq j$, such that
$$f_{\ell}^{t_\ell}\circ f_{\ell+1}^{t_{\ell+1}} \circ\dots \circ f_{j-1}^{t_{j-1}} \circ f_{j}^{t_j}(z_0) \in V_{\ell}-U_{\ell}$$
Let $\tilde{A}$ be a simple connected domain in $V_{i} - P_i$ which contains $\tilde{z_0} = f_{\ell}^{t_\ell}\circ f_{\ell+1}^{t_{\ell+1}} \circ\dots \circ f_{j-1}^{t_{j-1}} \circ f_{j}^{t_j}(z_0)$. Note that $f_{\ell}^{t_\ell}\circ f_{\ell+1}^{t_{\ell+1}} \circ\dots \circ f_{j-1}^{t_{j-1}} \circ f_{j}^{t_j}=f_{i}^a$, for some $a$, and so there exists  a simple connected domain $A$ such that $z_0 \in A$ and $f_{i}^{a}$ restricts to $A$ is an univalent map whose image is $\tilde{A}$, since $\func{f_{i}^{a}}{f_{i}^{-a}(V_{i}-P_{i})}{V_{i}-P_i}$ is a covering map. 

\begin{lem}\label{dist} In the conditions described above, the following holds:
\begin{enumerate}
\item Uniform expansion: We have:
$$\norm{ D(f_{i}^{a+1})(z_0))}_{\co - P} \geq \lambda^{j-i};$$

\item Distortion control: If $diam_{V_{i}-P_{i}}(\tilde{A}) \leq D$ then 
$$ \frac{1}{C(D)} \leq \frac{\norm{D(f_{i}^{a})(z_1)}_{\co - P}}{\norm{D(f_{i}^{a})(z_2)}_{\co - P}} \leq C(D),$$
for $z_1, z_2 \in A$.
\end{enumerate} \end{lem}
\begin{proof} The first statement is an immediate consequence of lemmas \ref{contra} and \ref{escontra}. For the second statement, note that $f_{j_\ell}^{t_\ell - 1}\circ f_{j_{\ell-1}}^{t_{\ell-1}} \dots \circ f_{j_1}^{t_1} \circ f_{j_0}^{t_0}=f_{i}^{a}$, for same $a$, and $f_{i}^a(z) \in f_{i}^{-1}(V_i - U_i)$. The map
$$f_{i}^{a}\colon f_{i}^{-a}(V_i - P_i) \rightarrow V_i - P_i$$ 
is a covering map, with $f_{i}^{-a}(V_i - P_i) \subset V_i - P_i$,  thus we can apply corollary \ref{disto}(twice) to obtain 2.
\end{proof} 

\begin{cor} The set $J(\mathbf{f})$ is  dense in $\co$.\end{cor}
\begin{proof} Let $z_0$ be a complex number which is not contained in $J(\mathbf{f})$. Then $z$ is not in $J(f_k)$ for sufficiently small $k$. Let $i$ be maximal such that  $z \in U_i$. For each small $k$ let $a(k)$ be minimal so that $f^{a(k)}(z_0) \in V_k - U_k$ and let $\gamma$ be the minimal geodesic between $f^{a(k)}_k(z_0)$ and $J(f_k)$ in $V_k - P_k$. The hyperbolic length of $\gamma$ is smaller than the constant $D$ in the lemma \ref{distcont}, since $f^{-1}_k(P_k) \subset J(f_k)$. By the previous lemma, the length of the  lift $f^{-a(k)}\gamma$ in the hyperbolic metric of $\co - P$ goes exponentially fast to zero, when $-k$ goes to infinity.
\end{proof}

\begin{cor}[rigidity] The towers in $\mathcal{T}(C_1,k,C_2)$ does not support non trivial invariant Beltrami fields.
\end{cor}
\begin{proof}  Let $\mu$ be an invariant line field to the tower $\mathbf{f}$. Because $K(f_i)$ does not support invariant line fields,  it is possible select a  point $z_0 \in \co - K(\mathbf{f})$ where $\mu$ is almost continuous. This means
$$ \lim_{\delta \rightarrow 0} \frac{\ell(\{z \colon |z - z_0| < \delta \text{ and } |\mu(z) - \mu(z_0)| \leq \epsilon   \})}{\ell(\{ z \colon |z - z_0| < \delta  \})}=1. $$
Here $\ell$ is the Lebesgue measure in $\co$. Since $J(\mathbf{f})$ is dense and by  small Julia sets everywhere theorem, for any $\alpha > 0$ there exists a polynomial like map $g_\alpha \colon V^\alpha_1 \rightarrow V^i_2$  so that:
\begin{itemize}
\item The map $g_\alpha \colon V^\alpha_1 \rightarrow V^\alpha_2$ belongs to $\mathcal{F}(f_i)$, for some $i$ (indeed, for any  $i$ small enough);
\item The map $g_\alpha \colon V^\alpha_1 \rightarrow V^\alpha_2$ belongs to $\mathcal{P}_n(C(C_1,C_2))$;
\item  $diam(J(g_\alpha)) \sim \alpha$;
\item $dist(z_0,J(g_\alpha)) \leq C(C_1,C_2) \cdot \alpha$.
\end{itemize}
Since $\mu$ is invariant by these maps, normalizing $g_\alpha$ so that $diam(J(g_\alpha))=1$, we can select a subsequence which converges to a polynomial like map which preserves a straight line field. This is absurd.  
\end{proof}

Let $\sigma = (\sigma_0,\sigma_1,\dots)$ be a sequence of m.c.d. We will denote by $R^{k}(\sigma)$ the bi-infinite sequence $$(\dots,\tilde{\sigma}_{-1},\tilde{\sigma}_0, \tilde{\sigma}_1,\dots)$$
where $\tilde{\sigma}_i = \sigma_{i+k}$ for $i \geq -k$. Fill the other positions in the sequence in an arbitrary way (we are interested in  convergent subsequences of $R^{k}(\sigma)$ when $k \rightarrow \infty$ in the space of bi-infinite sequences. Thus the other positions are not important for us). 

\begin{cor}\label{unique} There exists an unique bi-infinite tower  $\mathbf{g}_\sigma$ in $\mathcal{T}_n(C,1,C_2)$ with $C$-bounded combinatorics
 $$\sigma = (\dots, \sigma_{-2}, \sigma_{-1}, \sigma_{0}, \sigma_{1}, \sigma_2, \dots)$$
Here  unicity means that if $\mathbf{g}$ and $\mathbf{\tilde{g}}$ are bi-infinite towers with same combinatorics then the germs $g_i$ and $\tilde{g}_i$ are the same (up affine conjugacies). Furthermore there exists $C_1(C)$ such that the germ $g_i$ has a representation $g_i \colon U^i \rightarrow V^i$ which belongs to $\mathcal{P}_n^\infty(C,C_1)$. Notice that $C_1$ does not depend on $\sigma$.
\end{cor}
\begin{proof} 
\textbf{Existence:} Select a real infinitely renormalizable polynomial $p_0$ of type $n$ with combinatorics $\tilde{\sigma} = (\tilde{\sigma}_0, \tilde{\sigma}_1,\dots)$ so that for \textbf{any} $C$-bounded combinatorics there exists a sequence $k_i$ satisfying $R^{k_i}(\tilde{\sigma}) \rightarrow_{k} \sigma$. Using the complex bounds, select, for renormalizations deep enough, polynomial like representations in $\mathcal{P}_n^\infty(C,C(f))$. Then the finite tower $<g_{j,i}>_{|j|\leq k_i}$ defined by 
$$g_{j,i} = R^{k_i + j}p_0$$
has a subsequence which converges to a bi-infinite tower in $\mathcal{P}^\infty_n(C,C_1(f))$ with combinatorics $\sigma$.

\textbf{Unicity:} Let $\mathbf{f}$ and $\mathbf{g}$ be bi-infinity towers in $\mathcal{T}_n(C_1,1,C_2)$. Since $g_i\colon U_i^g \rightarrow V_i^g $ and $f_i \colon U_i^f \rightarrow V_i^f$ have the same combinatorics, there exists one $C(C_1,C_2)$-quasiconformal map  $\phi_i \colon \co  \rightarrow \co $ which maps $U_i^f$ in $U_i^g$ and it is a conjugacy  between $f_i$ and $g_i$ in $U_i^f$. When $i \rightarrow - \infty$ we have  $ V_i^f \rightarrow \co$. Thus $\phi_i$ admits a convergent subsequence to some quasiconformal map $\phi \colon \co \rightarrow \co$. This map is a conjugacy between the tower $\mathbf{f}$ and the tower $\mathbf{\tilde{g}}$, where $\tilde{g}_i$ is equal to $g_i$ restricts to $\phi(U_i^f)$.  Since the  Beltrami field $\frac{\overline{\partial} \phi}{\partial \phi}$ is invariant by the tower $\mathbf{f}$, the rigidity of towers implies that $\phi$ is conformal. Thus, up to affine maps, $\phi$ is the identity.\end{proof}

Let $B$ be a domain in $\co$. denote by $B(V)$ the Banach space of the holomophic functions defined in $V$ and with a continuous extension to $\partial V$. Denote by $\mathcal{A}^\infty_n(C)$ the set of germs $g$ in some level (and hence in the level $0$) of a bi-infinite tower in $\mathcal{T}_n(C,1,C_1)$, for some $C_1$.

\begin{thm}[convergence of renormalization]\label{conv} There exists $\delta = \delta(C_1)$ with the following property: For any $\epsilon > 0$ there exists $j_0 = j_0(C_1,C_2)$ so that if $f \in \mathcal{A}_n^\infty(C)$ and $g \ in \mathcal{P}_n^{\infty}(C_1,C_2)$ are polynomial like maps with same combinatorics then, for $j \geq j_0$:
\begin{itemize}
\item The germ $R^j(g)$ belongs to $\mathcal{P}_n^\infty(C_1,C_3)$. Here $C_3 = C_3(C_1)$.
\item The renormalizations $R^j(f)$, $R^j(g)$ belong to $B(\text{$\delta$-}K(R^j(f)))$, for $j \geq j_0$ and 
$$ |R^j(f)(z) - R^j(g)(z)| \leq \epsilon$$
for $z \in \delta \text{-}K(R^j(f))$.
\end{itemize}
\end{thm}
\begin{proof} Since $f \in \mathcal{P}_n^\infty(C_1,C(C_1))$ there exists a $K(C_2)$-quasiconformal conjugacy between  $f$ and $g$. Normalizing $R^j(f)$ and $R^{j}(g)$ so that $R^j(f)(0)=1=R^j(g)(0)$, we obtain quasiconformal conjugacies $\phi_j$ so that $\phi_j(0)=0$ and $\phi_j(1)=1$. We claim that $\phi_j$ converges uniformly in compact sets to identity. Indeed, suppose by contradiction there exist  sequences of maps $f_\ell$ and $g_\ell$, $f_\ell \in \mathcal{A}_n^\infty(C)$ and $g_\ell \in \mathcal{P}_n^{\infty}(C_1,C_2)$ , with same combinatorial type  so that the corresponding conjugacies $\phi_{j,\ell}$ does not converge in an uniform way to identity: in other words we can select  $\epsilon > 0$ so that $|\phi_{j_i, \ell_i}(z) - z| \geq \epsilon$, for some $z \in \co$ and with  $j_i \rightarrow \infty$. But lemma \ref{cbit} and proposition \ref{cconj} say that a subsequence of $\phi_{j_i, \ell_i}$ converges to a conjugacy between two bi-infinite towers, which do not support invariant line fields, so this conjugacy is a conformal map, hence it is the identity, which is a contradiction. 
Since $f \in \mathcal{P}_n^\infty(C_1,C(C_1))$, we can select representations  $R^{k}(f)\colon U^k \rightarrow V^k$ which belongs to $\mathcal{P}_n^\infty(C_1,\tilde{C}(C_1))$ and furthermore they are restrictions of iterates of $f$. This is possible for $k \geq k_0(C_1)$. Then $R^{k}(g)\colon \tilde{U}^k  \rightarrow  \tilde{V}^k$, where $\tilde{U}^k = \phi(U^k)$ and  $\tilde{V}^k = \phi(V^k)$, is a representation of $R^{k}(g)$. Since $2\delta$-$K(R^k(f)) \subset V^k$, for some $\delta = \delta(C_1)$, and $\phi_j$ is close to identity, one gets $\delta$-$K(R^k(g)) \subset \tilde{V}^k$ for $k \geq k_1(C_1,C_2)$. By lemmas \ref{comp1} and \ref{comp2}, $R^{k}(g) \in \mathcal{P}_n^\infty(C_1,C_3(C_1))$, which proves the first statement. To proof the second one, note that $\delta$-$K(R^k(f)) \subset \tilde{V}^k$ and $(\delta/2)$-$K(R^k(g)) \subset $ $\delta$-$K(R^k(f))$, for $k \geq k_2(C_2)$. By corollary \ref{delta} $R^{k}(g)\colon  R^{k}(g)^{-n}(\tilde{V}^k) \rightarrow R^{k}(g)^{-n + 1}(\tilde{V}^k)$ is a representation in $B(\delta \text{-}K(f))$, where $n = n(\delta/2,C_3)$. Since $\phi_j \rightarrow Id$ and the diameter of $\delta$-$K(R^k(f))$, after the normalization  $R^k(f)(0)=1$, is bounded yet, the proof is finished. 
\end{proof}

\subsection{Exponential convergence}\label{expo}  

\begin{defi} Let $\Lambda \subset \co$. We say that $z$ is a  $\delta$-\textbf{deep point} of $\Lambda$ if there exists $C$ so that  $B(\tilde{z},\tilde{r}) \in B(z,r)-\Lambda$ implies $\tilde{r} \leq C r^{1+\delta}$. Define $\mathcal{F}(f)$ by the set of functions $g\colon U \rightarrow V$ whose graph is contained in $\{(x,y) \colon f^i(x)=f^j(y)  \}$, for some $i,j \geq 0$.
\end{defi}

\begin{lem} The critical point  $0$ of $f \in P^{\infty}_n(C_1,C_2)$ is $\delta(C_1,C_2)$-deep.
\end{lem}
\begin{proof} Since the McMullen's proof  (proposition 8.8 in \cite{Mc2}), without modifications, works perfectly well in our situation, we will not repit it  here. \end{proof}

Let $\mathbb{H}^3$ be  the hyperbolic space and identify the Riemann sphere $\mathbb{S}$ with its ideal boundary. If $K$ is a subset of the Riemann sphere, denote by $hull(K)$ the set of points in $\mathbb{H}^3$ in geodesics which arrive in both directions in a point of $K$. Furthermore, given a quasiconformal vector field $v$ in $\mathbb{S}$, we said that $v$ is a \textbf{quasiconformal deformation} of a polynomial like map $f$ is $\overline{\partial} v$ is invariant by $f$. Define  the \textbf{visual distortion} $Mv$ by
$$ Mv(p) = \inf_{\overline{\partial} w=0} ||v - w||_\infty (p)    $$ 
The visual distortion measure the 'distance, of $v$ of the conformal vector fields. 
\begin{lem} For any $C_1, C_2$ there exists $r(C_1,C_2)$ with the following property: let  $f \in \mathcal{P}_n^{\infty}(C_1,C_2)$ and let $v$ be a quasiconformal deformation of $f$. Furthermore assume $S(p,r) \subset hull(K(f))$. Then
$$  Mv(p) \leq \frac{1}{2} sup_{q \in S(p,r)} Mv(q)  $$
\end{lem}
\begin{proof} The proof will be exactly as in lemma 9.12 in \cite{Mc2}, with some small modifications to avoid technical definitions: Suppose, by contradiction, there exist sequences  $r_i \rightarrow \infty$, $v_i$, $p_i \in \co$ and $f_i \in \mathcal{P}_n^{\infty}(C_1,C_2)$ so that 
\begin{itemize}
\item $S(p_i,r_i) \subset hull(K(f_i))$,
\item $Mv_i(p) \geq 1/2$,
\item $sup_{q \in S(p_i,r_i)} Mv_i(q) \leq 1 $.
\end{itemize}
We can assume that $p_i = p$. Then $S(p,r_i) \subset hull(K(f_i))$, with $r_i \rightarrow \infty$, which implies that, for all $z \in \overline{\co}$ and $\epsilon > 0$, $dist_{\overline{\co}}(z,K(f_i)) \leq \epsilon$, if $i$ is large enough. In particular, by small Julia sets everywhere lemma, for any $z \in \co$ and $d > 0$ there exists a sequence of polynomial-like maps $g_i \colon U_i \rightarrow V_i$, for $i$ large enough, so that 
\begin{itemize}
\item  $g_i \in \mathcal{F}(f_i)$,

\item  $g_i \in \mathcal{P}_n(C)$, $C = C(C_1,C_2)$,

\item $diam K(g_i) \sim d$,

\item $dist_{Eucl}(z,K(g_i)) \sim d$.
\end{itemize}
Because $sup_{q \in S(p,r_i)} Mv_i(q) \leq 1$, with $r_i \rightarrow \infty$, we can assume that the  sequence $v_i$, up to sums with conformal fields in the Riemann sphere, converge uniformly to a quasiconformal vector $v_\infty$. Moreover $\overline{\partial} v_i \rightarrow \overline{\partial} v_\infty$ as distributions. In particular $\mu_\infty = \overline{\partial} v_\infty$ is invariant by any limit of the sequence $g_i$. We claim that $\mu_\infty = 0$. Otherwise, select $z$ a point of almost continuity of $\mu_\infty$. Hence $\mu_\infty$ is almost a straight Beltrami field near to $z$, which is impossible since there are polynomial-like maps (which form a compact family after conjugacies by affine maps) in all scales so that $\mu_\infty$ is invariant for them. But this is a contradiction, since $Mv_i(p) \geq 1/2$ implies $Mv_\infty(p) \geq 1/2$, so $\mu_\infty \neq 0$. 
\end{proof}

\begin{thm}[Exponential convergence] Let $f \in \mathcal{A}_n^\infty(C_1)$ and let  $g$ be a map in  $\mathcal{P}_n^{\infty}(C_1,C_2)$ with the same combinatorics that $f$. There exist $k_0 = k_0(g)$ and $\delta$ so that $R^{k}(g) \in B(\delta-K(f))$, for $k \geq k_0$ and furthermore
$$ |R^{k}(f)(z) - R^{k}(g)(z)| \leq \alpha^{k}  $$
for $z \in $ $\delta$-$K(f)$ and $\alpha < 1$. Here $\delta$ and $\alpha$ depends only on $C_1$.
\end{thm}
\begin{proof}
 By the dynamic inflexibility theorem (\cite{Mc2}), if $\phi \colon \co \rightarrow \co$ is a $K$-quasiconformal map between $f\colon U_f \rightarrow V_f$ and $g\colon U_g \rightarrow V_g$, $f,g \in \mathcal{P}_n^{\infty}(C_1,C_2)$, with $\phi(U_f)=U_g$ then $\phi$ is $C^{1+\beta}$-conformal at $0$. It is not difficult to verify in the proof of dynamic inflexibility theorem that  
$$\beta = \beta( K, r(C_1,C_2), \delta(C_1,C_2))$$
Here $r$ is as in the previous lemma. Note that we can select $\phi$ such that $K=K(C_1,C_2)$. Hence $\phi$ satisfies
$$\phi(x) = \phi'(0) \cdot x + O(|x|^{1+\beta})$$
Since $diam K(R^k(f)) \leq C \cdot \lambda^k$ for some $C > 0$ and $\lambda < 1$, after normalize $R^k(f)$ so that $R^k(f)(0)=1=R^k(g)(0)$, $\phi$ define a conjugacy $\phi_k$ satisfying
$$ \phi_k(x) =   x + O(\alpha^{k})$$
for some $\alpha < 1$. Using arguments as in theorem \ref{conv}, the proof is finished. \end{proof}

\section{Apendice}

\subsection{A fixed point theorem} The following fixed point theorem was proved by de Melo and van Strien (\cite{MS}) when $K$ is a simplex.

\begin{pro}[de Melo-van Strien fixed point theorem \cite{MS}]\label{deMelo} Let $K$ be a bounded closed convex body in a finite dimensional normed linear space and let $T \colon \text{int }K \rightarrow \text{int }K$ be a continuous function such that
\begin{equation}
\lim_{x \rightarrow \partial K} \frac{|T(x) - x|}{dist(x,\partial K)} = \infty
\end{equation}
Then $T$ has a fixed point in $\text{int } K$.
\end{pro}

Let $K_1$ and $K_2$ be two bounded closed convex bodies in a finite dimensional normed linear space. The radial projection $\phi \colon K_1 \rightarrow K_2$ is defined by (1) $\phi(0)=(0)$ (2) If $x_i$, $i = 1,2$ are the unique points such that a ray beginning at $0$ crosses $\partial K_i$ then $\phi(\lambda x_1) = \lambda x_2$, for $\lambda > 0$. We need the following result

\begin{lem}[Sz.-Nagy-Klee Theorem: 
\cite{Kl}] The radial projection between $K_1$ and $K_2$ is bi-Lipschtizian.
\end{lem}

For an elementary proof (using gauge functions) and references, see \cite{Kl}.

\begin{figure}\label{fix}
\centering
\psfrag{0}{$0$}
\psfrag{d}{$a$}
\psfrag{b}{$b$}
\psfrag{a}{$\alpha(x)$}
\psfrag{x}{$x$}
\psfrag{T}{$\tilde{T}(x)$}
\includegraphics[width=0.30\textwidth]{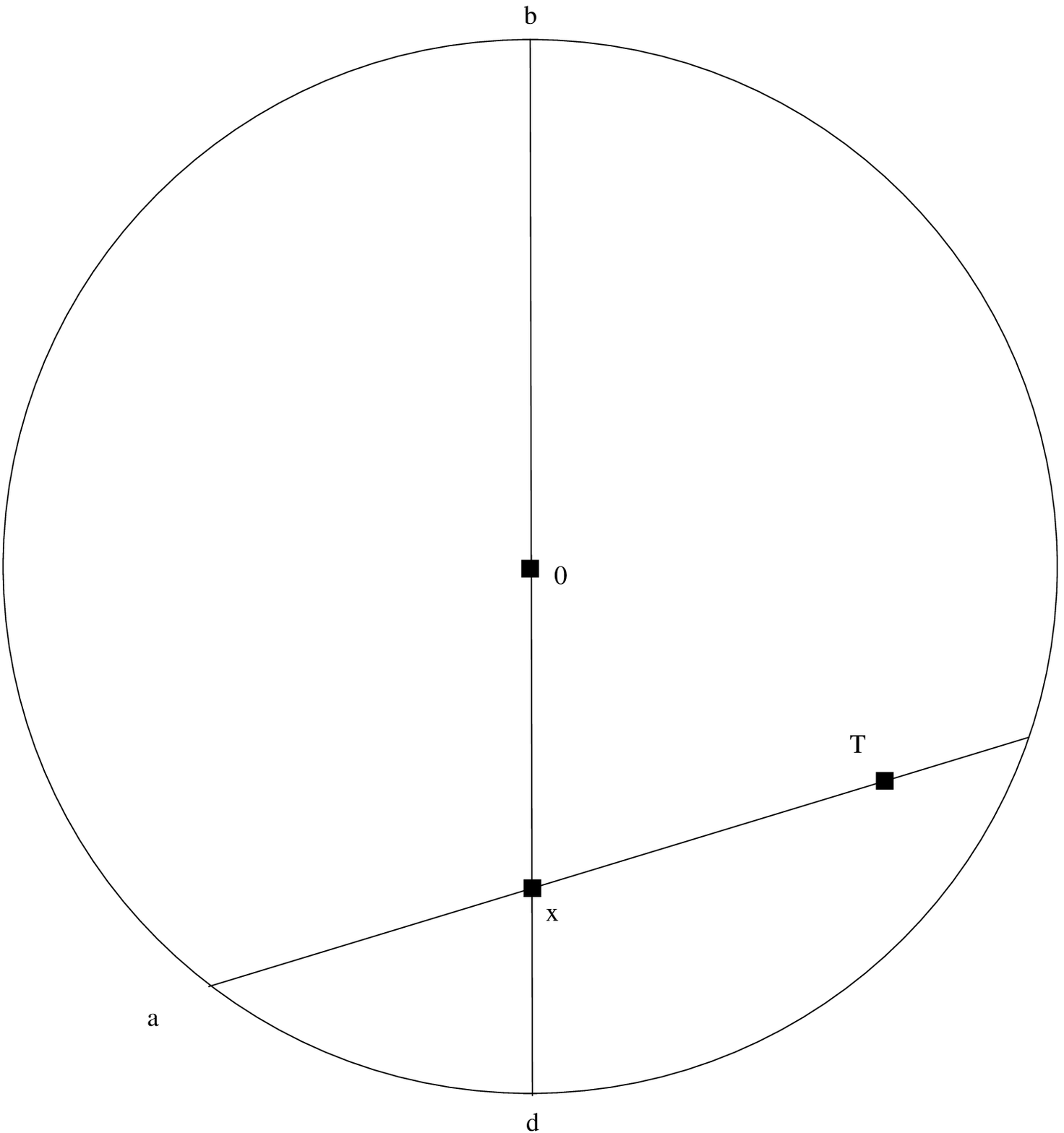}
\end{figure}

\begin{proof}[Proof of Proposition \ref{deMelo}] Suppose that $T$ does not have a fixed point in the interior of $K$. Let $\phi$ be the radial projection of unit ball $D = \{x \colon |x| \leq 1\}$ in $K$. Let $\tilde{T} = \phi^{-1} \circ T \circ \phi$. Then 
\begin{equation} 
\lim_{x \rightarrow \partial D} \frac{|\tilde{T}(x) - x|}{dist(x,\partial D)}= \infty 
\end{equation}
by the Sz.-Nagy-Klee Theorem. For a point in the interior of $D$, define the continuous function $\psi_\lambda(x) = (1-\lambda)x + \lambda \alpha (x)$, where $\alpha(x)$ is the unique point where the ray beginning at $\tilde{T}(x)$ and  containing $x$ crosses $\partial D$. We claim that $\psi_\lambda$ has a continuous extension to $D [0,1]$ such that $\phi_\lambda(x) = x$ for any $(\lambda, x) \in [0,1]  \partial D$, $\psi_0 = Id$ and $\psi_1(D) \subset \partial D$, which is absurd, since $\psi_\lambda$ will be a retraction of $D$ in its boundary. Indeed, if $a, b$  are as in the figure \ref{fix}, we have $|x - a|=dist(x,\partial D)$, $|x - b| = 1 -  dist(x,\partial D)$ and 
\begin{equation} 
(1 -  dist(x,\partial D))dist(x,\partial D) > |\tilde{T}(x) -x| |\alpha(x)  - x|
\end{equation}

Hence $|\psi_\lambda(x) - x|$ converges to zero when $x$ tends to $\partial D$ in an uniform way with respect to  $\lambda$. This proves the claim.
\end{proof}

\subsection{Hyperbolic domains on the plane} Let $\Omega$ be a hyperbolic domain on the plane and $\rho_\Omega|dz|$ its hyperbolic metric.  For $z \in \Omega$, define
$$  \beta_\Omega(z) = inf \{ | log \frac{|z-a|}{|b-a|}| \colon a, b \in \partial \Omega; |z-a|=dist(z,\partial \Omega) \} $$
To compare the Euclidean and hyperbolic metric on $\Omega$, we will use the following Beardon-Pommerenke results:

\begin{pro}[\cite{BePo}] There exists a constant $C$, with does not depend on $\Omega$, such that
$$  \frac{1}{2\sqrt{2}} \leq \rho_\Omega(z) \leq \frac{1}{dist(z,\partial \Omega)} \frac{C + \pi/4}{C +\beta_\Omega(z) }$$
\end{pro}

An annulus $A \subset \Omega$ is essential if the bounded component of $\mathbb{C} \setminus A$ contain point in $\partial \Omega$.

\begin{cor}[\cite{BePo}] The following holds:
\begin{itemize}
\item If $\Omega$ is a hyperbolic domain whose any essential annulus has modulus bounded by $M$, then there exists $C(M)$ such that
\begin{equation}\label{el} \frac{1}{C}  \frac{1}{dist(z,\partial \Omega)} \leq  \rho_\Omega(z) < C\frac{1}{dist(z,\partial \Omega)}
\end{equation}
For all  $z \in \Omega$;
\item If (\ref{el}) holds, then there exist $M(C)$ such that any essential annulus has modulus bounded by $M$.
\end{itemize}
\end{cor}

The domains satisfying the hypothesis of the previous corollary will be called \textbf{$M$-uniform domains}. Observe that if $\Omega$ is a uniform domain with maximum essential modulus bounded by $M$ and $D$ is a simply connect region in the plane, then $D \cap \Omega$ is also a uniform domain with the same bound for the maximal essential modulus.

\begin{pro}\label{ere} Let $\gamma$ be a Jordan curve  in a $M$-uniform domain $U$ with length $\ell \leq \ell_1$ in the hyperbolic metric of $U$ and let $D$ be the  bounded region in $\co - \gamma$. Then 
$$C_1(M,\ell_1) diam(D \cap \partial U) \leq dist(\gamma,D \cap \partial U) \leq C_2(M) diam(D \cap \partial U)   $$
Moreover
$$  diam(\gamma) \leq C_2(\ell_1, M) diam(D \cap \partial U) $$
\end{pro}
\begin{proof} Note that 
$$dist(\gamma,D \cap \partial U) \leq C_2(M) diam(D \cap \partial U)$$
otherwise there will be a large essential ring in $U$. Denote  $d = (1 + C_1)diam(D \cap \partial U)$ and fix $\lambda > 1$. Select an arbitrary $z \in D \cap \partial U$ and define 
$$ A_n = \{x \in \co \colon  \lambda^{n}d \leq dist(x,z)  \leq \lambda^{n+1}d   \}$$
Let $\alpha \subset U\cap A_n $, $n \geq 1$,  be a curve which touch both components of $\partial A_n$. Then the Euclidean length of $\alpha$ is at least $\lambda^n(\lambda -1)d$ and, if $\rho |dz|$ is the hyperbolic metric on $U$ then 
$$ \ell \geq  \int_0^{\lambda^n(\lambda -1)d} \rho(\gamma(t))|\gamma'(t)| dt \geq C(M) \int_0^{\lambda^n(\lambda -1)d} \frac{1}{dist(\gamma(t), D \cap \partial U)} dt \geq (1 - \frac{1}{\lambda})d$$ 
which proves the lemma. If the diameter of $\gamma$ is large relative to diameter of $D \cap \partial U$ then $\gamma$ crosses many rings $A_n$, so its hyperbolic length is large, which is absurd. To obtain the lower bound to  $dist(\gamma,D \cap \partial U)$,  notice that

$$  \gamma \subset \bigcup_{i \leq N} B(x_i,\frac{dist(x_i,\partial \Omega)}{2}) $$
for some $x_i$ in $\gamma$ and $N=N(\ell_1,M)$. It is easy to see that
$$diam(\gamma) \leq C(N) dist(x_i,\partial \Omega)$$
,for any i. Since $diam(\gamma) \geq diam(\partial \Omega \cap U)$, the proof is complete.
\end{proof}

\begin{rem} The previous lemma will be used in the following situation: Let $f \colon U \rightarrow V$ be an infinitely renormalizable polynomial-like map of type $n$ with bounded combinatorics. Then the postcritical set $P$ is a Cantor set with bounded geometry and hence $V - P$ is a $M$-uniform domain. Furthermore,  the hyperbolic length of the closed geodesics in $V - P$ is under control. Thus we can apply lemma \ref{ere}  for these geodesics.
\end{rem}

Choose $r < 1$ and large $m$ such that any annulus $R$ with $mod / R \geq m$, and any $z_0$ in the bounded component of $\co - R$, $\{z \in  \co \colon rC < |z-z_0| < C  \} \subset R$, for some $C$.

\begin{pro}[Distortion control] Let $\Omega_1$ and $\Omega_2$ be  $M$-uniform domains. There exists $C(M)$ with the following property: Let  $\func{g}{(D,z_0)}{(U,f(z_0))}$ be an univalent map such that $U \subset \Omega_2$ and  $D \subset \Omega_1$, where $D$ is a disc whose center is $z_0$. Then 
$$ \frac{1}{C} \leq \frac{\norm{g'(z_1)}_{\Omega_1, \Omega_2} }{ \norm{g'(z_2)}_{\Omega_1, \Omega_2}} \leq C $$
for any $z_1, z_2 \in A_m =\{z \in D \colon |z-z_0| \leq e^{-m} \cdot diam \ D/2   \}$. 
\end{pro}
\begin{proof} By the Koebe distortion theorem we immediately obtain
 $$\frac{1}{C_0(m)} \leq \frac{|g'(z_1)|}{|g'(z_2)|} \leq C_0(m)$$
Since 
$$\norma{g'(z)}{\Omega} = \frac{\rho_\Omega(g(z))}{\rho_\Omega(z) }|g'(z)|$$
it is suffice to prove that 
$$ \frac{1}{C_2(M,r)} \leq \frac{\rho_{\Omega}(z_1)}{\rho_{\Omega}(z_2)} \leq C_2(M,r) \text{ and } \frac{1}{C_3(M,r)} \leq  \frac{\rho_{\Omega}(g(z_1))}{\rho_{\Omega}(g(z_2))} \leq C_3(M,r) $$
Since $|z_i-z_0| \leq r dist(z,\partial \Omega)$, we gets
$$ \frac{1}{C_4(r)} \frac{1}{dist(z_0,\partial \Omega)} \leq \frac{1}{dist(z_i,\partial \Omega)} \leq  C_4(r) \frac{1}{dist(z_0,\partial \Omega)}$$ 
Note that $mod(g(A_r))=m$. Then $|f(z_i)-f(z_0)| \leq r dist(f(z_0),\partial \Omega)$ and we obtain
$$ \frac{1}{C_5(m)} \frac{1}{dist(g(z_0),\partial \Omega)} \leq \frac{1}{dist(g(z_i),\partial \Omega)} \leq  C_5(m) \frac{1}{dist(g(z_0),\partial \Omega)}$$
The last two equations and the $M$-uniformity of $\Omega$ proves the lemma.
\end{proof}

\begin{cor}\label{disto} Let $f\colon U \rightarrow \Omega_2$ be a covering map so that 
\begin{itemize}
\item  $\Omega_2$ is a  $M$-uniform domain;
\item  The domain $U$ is contained in a  $M$-uniform domain $\Omega_1$.
\end{itemize} 
If  $\tilde{A}$ is a simple connected domain inside $\Omega_2$ and $A$ is a connect component of $f^{-1}(\tilde{A})$, then 
$$ \frac{1}{C} \leq  \frac{|f'(z_1)|}{|f'(z_2)| } \leq C$$
and 
$$ \frac{1}{C} \leq  \frac{\norm{f'(z_1)}_{\Omega_1, \ \Omega_2}}{\norm{f'(z_2)}_{\Omega_1,\ \Omega2} } \leq C$$
for $z_1, z_2 \in A$ and $C= C(M,diam_{\Omega_2}(\tilde{A}))$.
\end{cor}
\begin{proof} Let $\gamma$ be the minimal geodesic in the hyperbolic domain $\Omega_2$ between $f(z_1)$ and $f(z_2)$. It is easy to see that there are $x_i \in \gamma$ such that $\gamma \subset \cup_{i \leq N} \{ z \in \co \colon |z-x_i| < r dist(x_i,\partial \Omega_2) \}$. Here $N=N(M,diam_{\Omega_2}(\tilde{A}))$. Thus we can apply the previous proposition (and the Koebe distortion lemma) in the inverse branches $g_i$ of $f$ in each ball $B_i = \{ z \in \co \colon |z-x_i| <  dist(x_i,\partial \Omega_2) \}$.  
\end{proof}

\subsection{Quasiconformal mappings} We say that a Jordan curve $J \subset \co$ is a \textbf{$\mathbf{C}$-quasicircle} if there is a $C$-quasiconformal map $\phi$ on the Riemann Sphere such that $\phi(S^1)=C$. 

\begin{lem} Let $\psi \colon \mathbb{D} \rightarrow \mathbb{D}$ be a $C_1$-quasiconformal map and $x, y \in \mathbb{D}$ with $dist_{\mathbb{D}}(x,y) \leq D$. Then there exists a $C_2(C_1,D)$-quasiconformal map $\tilde{\psi}\colon \mathbb{D} \rightarrow \mathbb{D}$ which coincide with $\psi$ in a neighborhood of $S^1$ and $\tilde{\psi}(x)=y$. 
\end{lem}
\begin{proof} It follows of lemma 5.2.3  in \cite{GW98} or the moving lemma at pg. 288 in \cite{Lyu}.
\end{proof}

\begin{pro}\label{quasi} Let $J$ be a $C_1$-quasicircle and $x$ a point in the bounded domain in $\co -J$ such that $dist(x,J) \geq \epsilon diam(J)$. Then there exists a $C_2(C_1,\delta)$-quasiconformal map \textbf{on the plane} $\tilde{\phi}$ such that $\tilde{\phi}(S^1)=J$ and $\tilde{\phi}(0)=x$.
\end{pro}
\begin{proof} Assume, without loss of generality, that $diam J = 1$. Consider a $C_1$- quasiconformal map in $\overline{\co}$ such that $\phi(S^1)=J$. After a composition with a Moebious transformation which preserves the circle, we can assume that $\phi(\infty)=\infty$. Furthermore, after translate and rotate $J$, we can assume that $\phi(0)=0$ and $\phi(1)=1$. Since the set of $C_1$-quasiconformal maps on the plane such that  $\phi(0)=0$ and $\phi(1)=1$ is compact, there exists $\delta > 0$ such that for $a, b \in \overline{D}$, $|a-b| \leq \delta$ implies $|\phi(a) - \phi(b)| \leq \epsilon$. In particular $dist(\phi^{-1}(x), S^1) \geq \delta$. By the previous lemma, we obtain a $C_2(C, \delta)$-quasiconformal map $\tilde{\phi}$ on the plane which is equal to $\phi$ outside $\mathbb{D}$ and $\tilde{\phi}(0)=x$.
\end{proof}

\section{Acknowledge} 
I wish to thank W. de Melo for introducing me on this subject and for useful conversations about the mathematics and the style of this article. I also grateful to IMPA, where this work was done.

\end{document}